\documentclass[12pt,english]{amsart}
\usepackage[T1]{fontenc}
\usepackage[latin1]{inputenc}
\usepackage{geometry}
\usepackage{amssymb}

\makeatletter
 \theoremstyle{plain}    
 \newtheorem{thm}{Theorem}[section]
 \numberwithin{equation}{section} 
 \numberwithin{figure}{section} 
 \theoremstyle{plain}
 \theoremstyle{plain}    
 \newtheorem{cor}[thm]{Corollary} 
 \theoremstyle{definition}
  \newtheorem{example}[thm]{Example}
 \theoremstyle{remark}
 \newtheorem{rem}[thm]{Remark}
 \theoremstyle{plain}    
 \newtheorem{prop}[thm]{Proposition} 
 \theoremstyle{plain}    
 \newtheorem{lem}[thm]{Lemma} 
 \theoremstyle{remark}    
 \newtheorem{claim}[thm]{Claim}
 \theoremstyle{plain}    
 \newtheorem{conjecture}[thm]{Conjecture} 

\def\makebbb#1{
    \expandafter\gdef\csname#1\endcsname{
        \ensuremath{\Bbb{#1}}}
}
\makebbb{R}
\makebbb{N}
\makebbb{Z}
\makebbb{C}
\makebbb{G}
\makebbb{E}
\makebbb{H}
\makebbb{P}
\makebbb{B}
\makebbb{K}
\makebbb{Q}

\usepackage{babel}
\makeatother
\begin{document}

\title{Bergman kernels and equilibrium measures for polarized pseudoconcave
domains}

\author{Robert Berman}

\curraddr{Institut Fourier, 100 rue des Maths, BP 74, 38402 St Martin d'Heres
(France)}

\email{robertb@math.chalmers.se}

\begin{abstract}
Let $X$ be a domain in a closed polarized complex manifold $(Y,L),$
where $L$ is a (semi-)positive line bundle over $Y.$ Any given Hermitian
metric on $L$ induces by restriction to $X$ a Hilbert space structure
on the space of global holomorphic sections on $Y$ with values in
the $k$th tensor power of $L$ (also using a volume form $\omega_{n}$
on $X).$ In this paper the leading large $k$ asymptotics for the
corresponding Bergman kernels and metrics are obtained in the case
when $X$ is a pseudoconcave domain with smooth boundary (under a
certain compatibility assumption). The asymptotics are expressed in
terms of the curvature of $L$ and of the boundary of $X.$ The convergence
of the Bergman metrics is obtained in a very general setting where
$(X,\omega_{n})$ is replaced by any measure satisfying a Berstein-Markov
property. As an application the (generalized) equilibrium measure
of the polarized pseudoconcave domain $X$ is computed explicitely.
Other applications to the zero and mass distribution of random holomorphic
sections and the eigenvalue distribution of Toeplitz operators will
appear elsewhere.
\end{abstract}


\maketitle
\tableofcontents{}

\section{Introduction}

Let $L$ be a holomorphic line bundle over a closed (i.e. compact
without boundary) projective complex manifold $Y$ of dimension $n.$
Denote by $H^{0}(Y,L^{k})$ the vector space of all global holomorphic
sections on $Y$ with values in the $k$th tensor power of $L.$ Any
given Hermitian metric $\phi$ on $L$ and a domain $X$ in $Y$ together
with a volume form $\omega_{n}$ induces an $L^{2}-$norm on $H^{0}(Y,L^{k})$
obtained by integrating the point-wise norms of sections in $H^{0}(X,L^{k})$
over the domain $X.$ The corresponding Hilbert space will be denoted
by $\mathcal{H}(Y,L^{k})_{X}.$ The \emph{Bergman kernel} $K^{k}(x,y)$
of the Hilbert space $\mathcal{H}(Y,L^{k})_{X}$ is the integral kernel
of the orthogonal projection from the space of all smooth sections
with values in $L^{k}$ onto $\mathcal{H}(Y,L^{k})_{X}.$

In this paper the situation when the curvature form $dd^{c}\phi$
is \emph{semi-positive} and the domain $X=\{\rho\leq0\}$ has a smooth
strictly pseudo-concave boundary, i.e. the Levi curvature form $dd^{c}\rho$
of the boundary is \emph{negative,} will be mainly investigated. Then
$X$ (or rather the triple $(X,L,\phi$)) will be called a \emph{polarized
pseudo-concave domain.}

In the case when $X=Y$ and the curvature form $dd^{c}\phi$ is \emph{positive}
the asymptotics of the Bergman kernel $K_{k}(x,y)$ as $k$ tends
to infinity have been studied extensively \cite{ti,z,d-l-m,b-b-s}
and are by now very well-understood due to strong locazation properties.
For example, in scaled coordinates on {}``length-scales'' of the
order $1/k^{1/2}$ the Bergman kernels $K_{k}(x,y)$ converge (with
all derivatives) to constant curvature model kernels. In particular,
the leading asymptotics of the \emph{Bergman measure} $B^{k}\omega_{n},$
where $B^{k}(y):=\left|K^{k}(y,y)\right|_{k\phi}^{2}$ (the point-wise
norm) may be expressed in terms of the local curvature of $L:$ \begin{equation}
k^{-n}B^{k}\omega_{n}\rightarrow(dd^{c}\phi)^{n}/n!\label{eq:intro bergman measure}\end{equation}
uniformly on $Y.$ As an immediate consequence Tian's almost isometry
theorem \cite{ti} holds \begin{equation}
k^{-1}\Omega_{k}:=k^{-1}dd^{c}\textrm{ln\,$K^{k}(y,y)$}\rightarrow dd^{c}\phi\label{eq:intro berm metric conv}\end{equation}
uniformly on $Y,$ where $k^{-1}\Omega_{k}$ is called the (normalized)
$k$ th \emph{Bergman metric on $Y.$} Note that the latter asymptotics
are considerably weaker than \ref{eq:intro bergman measure}.

One notable application of these asymptotics was introduced by Shiffman-Zelditch
in their study of random zeroes of random and quantum chaotic holomorphic
sections \cite{sz2} (see section \ref{sub:Relations-to-random} below)
and was further developed in a series of papers (for example with
Bleher \cite{b-s-z,b-s-z2}).

A concrete realization of the situation studied in this paper is obtained
by taking $Y$ as the $n-$dimensional projective space $\P^{n}$
and $L$ as the hyperplane line bundle $\mathcal{O}(1)$. Then the
Hilbert space $\mathcal{H}(Y,L^{k})_{X}$ may be identified with the
space of all polynomials $p_{k}(z)$ in $\C^{n}$ of total degree
at most $k,$ equipped with the weighted norm \begin{equation}
\left\Vert p_{k}\right\Vert _{k\phi,X}^{2}:=\int_{X}\left|p_{k}(z)\right|^{2}e^{-k\phi(z)}\omega_{n},\label{eq:intro norm}\end{equation}
where $X$ has been replaced by its restrictio to the affine piece
$\C^{n},$ where $\phi$ is a smooth plurisubharmonic function of
logarithmic growth and $\omega_{n}$ is the restricted Fubini-Study
volume form (then the integrals are finite). Moreover, $X$ is by
assumption the complement of a bounded pseudoconvex domain in $\C^{n}.$
Such {}``weighted polynomials'' (with $\omega_{n}$ replaced by
a measure supported on a {}``arbitrary'' set $X)$ have been recently
studied in various contexts. See for example \cite{b-l2} and Bloom's
appendix in the book \cite{sa-to} by Saff-Totik and the book \cite{dei}
by Deift for the case when $E$ is a set in $\E,$ concerning relations
to (hermitian) random matrix theory. Very recently Bloom-Shiffman
\cite{bl-sh} studied the {}``unweighted'' situation obtained by
setting $\phi=0$ in \ref{eq:intro norm} and letting $X$ be a {}``regular''
bounded set in $\C^{n}$ (see section \ref{sub:Regularity}). Using
pluripotential theory \cite{kl} it was shown in \cite{bl-sh} that
the corresponding normalized $k$ th {}``Bergman volume form'' (compare
formula \ref{eq:intro berm metric conv}) converges weakly to the
\emph{equilibrium measure} $\mu_{e}$ of $X,$ supported on the (Shilov-)
boundary of $X:$ \begin{equation}
(\Omega_{k}/k)^{n}/n!\rightarrow\mu_{e}\label{eq:intro s-b}\end{equation}
When the domain $X$ is polarized (i.e. $dd^{c}\phi>0$) the situation
in the \emph{interiour} of $X$ can be shown to localize (as in \ref{eq:intro bergman measure}).
The main purpose of the present paper is to study the influence of
the \emph{boundary} on the Bergman kernel asymptotics of $\mathcal{H}(Y,L^{k})_{X}$
and on a generalized equilibrium measure of the polarized pseudoconcave
domain $X$ (defined following the very recent work \cite{g-z} of
Guedj-Zeriahi). In the situation of Shiffman-Bloom refered to above
these objets may, in general, not be expressed in terms of the local
curvature of the boundary $\partial X$ of the domain $X.$ However,
under the assumption of global negativitity of the curvature of the
boundary $\partial X$ there is a natural locally defined candidate
for the boundary contribution, namely the following $2n-1$ form,
invariantly defined on the boundary of $X:$ 

\begin{equation}
\mu:=\int_{0}^{T}(dd^{c}\phi+tdd^{c}\rho)^{n-1}\wedge d^{c}\rho)dt/(n-1)!,\label{eq:def of mu}\end{equation}
 where $T$ is the following function on $\partial X,$ that will
be referred to as the \emph{slope function}: \[
T=\sup\left\{ t\geq0:\,(dd^{c}\phi+tdd^{c}\rho)_{x}\geq0\,\textrm{along }T^{1,0}(\partial X)_{x}\right\} .\]
The point is that $T$ is finite when $\partial X$ is pseudoconcave.
It will be shown that, further assuming a certain \emph{compatibility}
between the curvature $dd^{c}\rho$ of the boundary $\partial X$
and the curvature $dd^{c}\phi$ of line bundle $L,$ leads to localization
properties of the Bergman kernel asymptotics and the (generalized)
equilibrium measure. In fact, as illustrated by the examples in section
\ref{sub:Counter-examples}, there are large classes of polarized
pseudo-concave domains $X$ where the localization properties hold
precisely when the assumption on {}``compatible curvatures'' holds. 

The main results below are based on the Bergman kernel asymptotics
obtained in section \ref{sec:Bergman-kernel-asymptotics}. A major
role in the proofs of these asymptotics is played by the local holomorphic
Morse inequalities obtained in \cite{berm1,berm3}. In the present
setting these inequalities can be seen as refined versions of the
Bernstein-Markov inequalities used by Shiffman-Bloom (compare section
\ref{sub:Bernstein-Markov-measures-and}). In the last section some
open problems concerning general smooth domains $X$ (and even more
general situations) are formulated. These open problems should be
seen in the light of some very recent developments that have appeared
since the preprint of the first version of the present paper appeared:
in \cite{berm4b} the situation when $X=Y,$ but the curvature of
$L$ is arbitrary is studied and in \cite{b-l2,b-l3} the planar case
is studied. 

Finally we turn to the precise statement of the main results (see
section \ref{sec:Setup} for further notation).

\subsection{Overview of the main present results }

The polarized pseudoconcave domain $X$ will be said to have {}``compatible
curvatures'' when the slope function $T$ above is constant on $\partial X$
for some choice of the defining function $\rho$ and certain further
assumptions depending on the {}``filling'' $Y-X$ of $X$ hold (see
section \ref{sub:Assumptions-on-the}). For example, in the case of
polynomials refered to above the compatibility assumption holds if
$-\rho=\phi$ in \ref{eq:intro norm}.

\subsubsection*{Bergman kernel asymptotics (section \ref{sec:Bergman-kernel-asymptotics})}

The first main result gives the convergence as a measure of the Bergman
kernel:

\begin{thm}
\label{thm:conv of K to current}Let $K^{k}$ be the Bergman kernel
for the Hilbert space $\mathcal{H}(Y,L^{k})_{X}$ associated to the
polarized pseudoconcave domain $X$ with compatible curvatures. Denote
by $\Delta_{X\times X}$ and $\Delta_{\partial X\times\partial X}$
the currents of integration on the diagonal in $X\times X$ and $\partial X\times\partial X,$
respectively. Then the sequence of measures \[
\begin{array}{lr}
k^{-n}\left|K^{k}(x,y)\right|_{k\phi}^{2}1_{X}(x)\omega_{n}(x)\wedge1_{X}(y)\omega_{n}(y)\end{array}\]
converges on $Y\times Y$ to \[
[\Delta_{X\times X}]\wedge1_{X(0)}(dd^{c}\phi)^{n}/n!+[\Delta_{\partial X\times\partial X}]\wedge\mu\]
in the weak {*}-topology, where $\mu$ is the $2n-1$ form \ref{eq:def of mu}
on $\partial X.$
\end{thm}
In fact, in order to prove the previous theorem the following {}``special
case'' will first be shown for the corresponding Bergman measure
(compare formula \ref{eq:intro bergm measure}): \begin{equation}
k^{-n}B^{k}1_{X}\omega_{n}\rightarrow1_{X}(dd^{c}\phi)_{n}+[\partial X]\wedge\mu\label{eq:intro bergm measure}\end{equation}
 weakly as measures on $Y.$ The next theorem concerns the scaling
convergence of the Bergman kernel $K^{k}$ close to the diagonal.
It shows that after scaling $K^{k}$ converges to constant curvature
model kernels (at least after choosing a subsequence). The scalings
are expressed in terms of the {}``normal'' local coordinates introduced
in section \ref{sub:int reg} and \ref{sub:bd reg}, respectively.
In the statement below the dependence on the fixed center (which is
the point $x$ in the interior and the point $\sigma$ at the boundary)
has been suppressed. 

\begin{thm}
\label{thm:scal conv for K}Let $K^{k}$ be the Bergman kernel for
the Hilbert space $\mathcal{H}(Y,L^{k})_{X}$ associated to the polarized
pseudoconcave domain $X$ with compatible curvatures. $K^{k}$ has
a subsequence $K^{k_{j}}$ such that for almost any point $x$ in
the interior of $X$ (i.e. $x\in X-E,$ where $E$ has measure zero)
the following scaling asymptotics hold in the $\mathcal{C}^{\infty}-$topology
on any compact subset of $\C_{z}^{n}\times\C_{z'}^{n}:$\[
(i)\, k_{j}^{-n}K^{k_{j}}(z/k_{j}^{1/2};z'/k_{j}^{1/2})\rightarrow K^{0}(z;z'),\]
 where $K^{0}$ is the corresponding model Bergman kernel (formula
\ref{eq:model K int}). Moreover, for almost any fixed point $\sigma$
in the boundary $\partial X$ (i.e. $\sigma\in\partial X-F,$ where
$F$ has measure zero in $\partial X$) the following scaling asymptotics
hold in the $\mathcal{C}^{\infty}-$topology on any compact subset
of $\C_{z,w}^{n}\times\C_{z',w'}^{n}:$ \[
(ii)\, k_{j}^{-(n+1)}K^{k_{j}}(z/k_{j}^{1/2},w/k_{j};z'/k_{j}^{1/2},w'/k_{j})\rightarrow K^{0}(z,w;z',w'),\]
where $K^{0}$ is the corresponding model Bergman kernel (formula
\ref{eq:model K bd}). Furthermore, the same statement holds after
replacing $K^{k}$ with any subsequence $K^{k_{l}}$ (a priori $E$
and $F$ then depend on the subsequence $K^{k_{l}}).$ 
\end{thm}
The model kernel $K^{0}$ associated to a point in the boundary may
be expreesed by the following suggestive formula, where $\rho_{0}$
denotes the (polarized) defining function of the corresponding constant
curvature model domain:\[
K^{0}=\frac{1}{4\pi}\frac{1}{\pi}\textrm{det}(dd^{c}\rho_{0})e^{\phi_{0}}P(\frac{\partial}{\partial\rho_{0}})\frac{\partial}{\partial\rho_{0}}(\frac{e^{T\rho_{0}}-1}{\rho_{0}}),\]
 where $P$ is the characteristic polynomial of the linear operator
$\{ dd^{c}\phi\}_{x}\{-dd^{c}\rho\}_{x}^{-1}.$ This kernel should
be compared with the one obtained by Shiffman-Zelditch \cite{sz3}
in the \emph{one-}dimensional unweighted case refered to above (the
later kernel is essentially given by $\frac{e^{v}-1}{v}$ in special
coordinates). The proofs in \cite{sz3} relied on classical results
of Carleman concerning the corresponding orthogonal polynomials and
the exterior Riemann mapping theorem. The corresponding unweighted
higher-dimensional scaling result in $\C^{n}$ was stated as an open
problem in \cite{bl-sh}.

\subsubsection*{Bergman metric asymptotics (section \ref{sec:Bergman-metric-asymptotics})}

Denote by $F_{k}$ the interior scaling maps on $\C^{n}$, as well
as the boundary ones, corresponding to the scaling of the coordinates
in theorem \ref{thm:scal conv for K} above. The following theorem
gives the convergence of the $k$th Bergman metric on $Y$ induced
by the polarized pseudoconcave domain $X$ (compare section \ref{sec:Bergman-metric-asymptotics}
for definitions). 

\begin{thm}
\label{thm:Berg metric}Let $\Omega_{k}$ \emph{be the Bergman metric
on $Y$} \emph{induced by the polarized domain $X$ with compatible
curvatures.} \emph{Then the following convergence holds for the corresponding
normalized volume form:} \[
(\Omega_{k}/k)_{n}\rightarrow1_{X}(dd^{c}\phi)_{n}/n!+[\partial X]\wedge\mu,\]
when $k$ tends to infinity, as measures in the weak{*}-topology,
where $\mu$ is the $2n-1$ form \ref{eq:def of mu} on $\partial X.$

Moreover, the following scaling asymptotics for the $p$th exterior
power of $\Omega_{k}$ hold (after replacing $K_{k}$ with a subsequence
as in theorem \ref{thm:scal conv for K}) around almost any interior
point: \[
(i)\, F_{k}^{*}(\Omega_{k})^{p}\rightarrow(dd^{c}\phi)^{p}\]
(with uniform convergence on each compact set) and around almost any
boundary point: \[
(ii)\, F_{k}^{*}(\Omega_{k})^{p}\rightarrow(dd^{c}\phi+tdd^{c}\rho+dt\wedge d^{c}\rho)^{p}\]
(with uniform convergence on each compact set), where $t$ is the
following function of $\rho:$ $t=\frac{\partial}{\partial\rho}\textrm{ln}B^{0}(\rho)$
(see formula \ref{eq:model B bd}) so that $t(-\infty)=0$ and $t(\infty)=T$
(where $T$ is the slope function in formula \ref{eq:def of mu}).
\end{thm}

\subsubsection*{Equilibrium measures (section \ref{sec:Distribution-of-random})}

Following the recent work \cite{g-z} of Guedj-Zeriahi (see also \cite{b-l2}
for the weighted case in $\C^{n})$ let now $X$ be any compact set
in $X$ and $\phi$ the {}``restriction'' to $X$ of a continuous
metric on $L.$ The corresponding \emph{equilibrium metric} on $L\rightarrow Y$
is defined by

\begin{equation}
\phi_{e}(y)=\sup\left\{ \widetilde{\phi}(y):\,\widetilde{\phi}\in\mathcal{L}_{(X,L)},\,\widetilde{\phi}\leq\phi\,\,\textrm{on$\, X$}\right\} .\label{eq:extem metric}\end{equation}
 where $\mathcal{L}_{(X,L)}$ is the class consisting of all (possibly
singular) metrics on $L$ with positive curvature current. Consider
the {}``regular'' case when $\phi_{e}$ is in $\mathcal{L}_{(X,L)}$
(compare section \ref{sub:Regularity}). The Monge-Ampere measure
$(dd^{c}\phi_{e})^{n}/n!$ is called the \emph{equilibrium measure}
associated to $(X,\phi).$ It was recently introduced in the more
general global setting of quasiplurisubharmonic functions by Guedj-Zeriahi
\cite{g-z}, building on the work of Bedford-Taylor, Demailly and
others. The item $(i)$ in the following theorem implies that if $Y$
is \emph{any} smooth domain then the normalized $k$ th Bergman volume
form converges to the corresponding equilibrium measure (see section
\ref{sub:Bernstein-Markov-measures-and} for the definition of Bernstein-Markov
measures etc). In item $(ii)$ the optimal rate of convergence (saturated
by the model examples in section \ref{sec:Model-examples} - see \cite{berm3})
is obtained in the case when $Y$ is strongly pseudoconcave.

\begin{thm}
\label{thm:ln k is equil}Let $X$ be a compact set in $Y$ and $\phi$
a continous metric on an ample line bundle $L\rightarrow Y.$ 

$(i)$ Let $\omega_{n}$ be a volume form on $Y.$ If $1_{X}\omega_{n}$
has the Bernstein-Markov property w.r.t $(X,\phi),$ then the following
uniform convergence holds on all of $Y:$ \begin{equation}
k^{-1}\textrm{ln\,$K^{k}(y,y)\rightarrow\phi_{e}(y)$}\label{eq:statement theorem ln k}\end{equation}
 where $K^{k}$ is the Bergman kernel associated to $(X,\omega_{n},\phi).$
In particular, the equilibrium metric $\phi_{e}$ is continuous then,
i.e. $(X,\phi)$ is regular then.

$(ii)$ If furthermore $X$ is assumed to be a pseudoconcave domain
with smooth boundary and $\phi$ is smooth, then the rate of the convergence
in \ref{eq:statement theorem ln k} is of the order $(n+1)\ln k/k.$ 

$(iii)$ If $\nu$ is any fixed measure which has the Bernstein-Markov
property w.r.t $(X,\phi)$ and $(X,\phi)$ is regular, then the uniform
convergence \ref{eq:statement theorem ln k} holds for the Bergman
kernel associated to $(X,\nu,\phi).$ 

Moreover, if $L$ is only assumed to be a semi-positive line bundle,
then the  following convergence holds under any of the assumptions
$(i),(ii)$ or $(iii)$ above for the normalized volume form of the
corresponding $k$th Bergman metric \emph{$\Omega_{k}:$\begin{equation}
(\Omega_{k}/k)^{n}\rightarrow(dd^{c}\phi_{e})^{n}\label{eq:bergman vol as equi}\end{equation}
}when $k$ tends to infinity, as measures in the weak{*}-topology. 
\end{thm}
The theorem above generalizes the result \ref{eq:intro s-b} of Shiffman-Bloom,
concerning the unweighted case in $\C^{n}$ (as well as Theorem 2.1
in \cite{bl2} concerning the weighted case for $X$ a compact set
in $\C^{n}).$ The proof is similar to Demailly's $\overline{\partial}-$proof
of Siciak's fundamental convergence result for the $L^{\infty}-$
version of the Bergman metrics (compare remark \ref{rem:l infty})
in the unweighted case in $\C^{n}$ \cite{d3}. See also \cite{g-z}
for the global \emph{polarized} case of this $L^{\infty}-$ version
of the result. Also note that in the case when $X=Y$ the convergence
towards the equilibrium measure was obtained in \cite{berm4b} for
\emph{any} line bundle $L$ (when suitably formulated). The proof
of the lower bound in the convergence of the theorem above uses the
Ohsawa-Takegoshi extension theorem, which allows a precise controle
on the rate of the convergence. 

In case $X$ is a polarized pseudoconcave domain that satisfies the
assumption about compatible curvatures (section \ref{sub:Assumptions-on-the})
the equilibrium measure can now be computed explicitly using theorem
\ref{thm:Berg metric} (without assuming that $L$ is ample): 

\begin{cor}
\label{cor:expl equil}Let $X$ be a polarized pseudoconcave domain
with compatible curvatures (section \ref{sub:Notation-and-setup}).
\emph{Then the (generalized) equilibrium measure $(dd^{c}\phi_{e})^{n}/n!$
of the polarized domain $X$ is given by} \[
(dd^{c}\phi_{e})/n!=1_{X}(dd^{c}\phi)^{n}/n!+[\partial X]\wedge\mu,\]
where $\mu$ is the $2n-1$ form \ref{eq:def of mu} on $\partial X.$
\end{cor}

\subsection{\label{sub:Relations-to-random}Relations to random sections and
Toeplitz operators}

In a sequel \cite{berm5} to his paper the present results will be
applied to the study of various random measure processes. The starting
point is that any Hilbert space $\mathcal{H}_{k}$ (here $\mathcal{H}(Y,L^{k})_{X}$)
comes equipped with a natural Gaussian probability measure. As shown
by Shiffman-Zelditch the Berfman measure $k^{-n}B^{k}\omega_{n}$
(formula \ref{eq:intro bergm measure}) then represents the \emph{expected
mass distribution} $\E(\left|f_{k}\right|^{2}\omega_{n})$ of a random
section $f_{k}$ in $\mathcal{H}_{k}$ and the Bergman volume form
$(dd^{c}(\textrm{ln\,$K^{k}(z,z)))^{n}/n!$}$ represents the \emph{expexted
distribution of simultaneous zeroes of $n$ random sections} in $\mathcal{H}_{k}.$
Moreover, the \emph{}variance of the mass distribution can to the
leader order be expressed in terms of the eigenvalue distribution
of \emph{}Toeplitz operators \emph{}acting on $\mathcal{H}_{k}$ (compare
\cite{sz2}), which in turn may be obtained from the weak convergence
of the measure $k^{-n}\left|K^{k}(x,y)\right|_{k\phi}^{2}\omega_{n}(x)\wedge\omega_{n}(y)$.
Furthermore, the scaling properties of the Bergman kernel $K^{k}(x,y)$
are used to express the limit correlations between random zeroes (compare
\cite{b-s-z,b-s-z2}). In \cite{berm5} the non-local effects appearing
when the condition about {}``compatible curvatures'' does not hold
will also be investigated and related to the situation studied in
\cite{berm4b}, as well as the Hele-Shaw flow in interface dynamics
(also called Laplacian growth) \cite{za}.

\section{\label{sec:Setup}Setup}

\subsection{\label{sub:Notation-and-setup}Notation}

Let $L$ be an Hermitian holomorphic line bundle over a complex manifold
$Y$. The Hermitian fiber metric on $L$ will be denoted by $\phi.$
In practice, $\phi$ is considered as a collection of local functions.
Namely, let $s$ be a local holomorphic trivializing section of $L,$
then locally, $\left|s(z)\right|_{\phi}^{2}=e^{-\phi(z)}.$ If $\alpha_{k}$
is a holomorphic section with values in $L^{k},$ then it may be locally
written as $\alpha_{k}=f_{k}s^{\otimes k},$ where $f_{k}$ is a local
holomorphic function and the point-wise norm of $\alpha_{k}$ may
be written as\begin{equation}
\left|\alpha_{k}\right|_{k\phi}^{2}=\left|f_{k}\right|^{2}e^{-k\phi(z)}.\label{eq:ptwise norm}\end{equation}
 The canonical curvature two-form of $L$ can be globally expressed
as $\partial\overline{\partial}\phi$ and the normalized curvature
form $i\partial\overline{\partial}\phi/2\pi=dd^{c}\phi$ (where $d^{c}:=i(-\partial+\overline{\partial})/4\pi)$
represents the first Chern class $c_{1}(L)$ of $L$ in the second
real de Rham cohomology group of $X$ \cite{gri}. A line bundle will
be said to be \emph{(semi-) positive} if there is some smooth metric
$\phi$ on $L$ with (semi-) positive curvature form (i.e. the matrix
$(\frac{\partial^{2}\phi}{\partial z_{i}\partial\bar{z_{j}}})$ is
(semi-) positive). 

Let $X$ be a smooth strictly pseudoconcave domain in $Y.$ This means
that there is a defining function $\rho$ (i.e. $X=\left\{ \rho\leq0\right\} $
and $d\rho\neq0$ on $\partial X$) such that the restriction of the
Levi curvature form $\partial\overline{\partial}\rho$ to the maximal
complex subbundle $T^{1,0}(\partial X)_{x}$ of the real tangentbundle
of $\partial X$ is negative (i.e. the Levi curvature of $\partial X$
is negative). The degenerate case $X=Y$ is allowed in the previous
definition of $X$ and corresponds to the situation studied in \cite{sz2,b-s-z,b-s-z2}
(when $dd^{c}\phi$ is \emph{strictly} positive).

We will assume that $Y$ is a projective manifold with a semi-positive
line bundle $L$ (which is positive at some point in $Y).$ In case
the curvature is positive on all of $Y,$ the pair $(Y,L)$ it usually
called a \emph{polarized manifold} in the literature. Fix a (possibly
singular) Hermitian metric $\phi$ on $L$ over $Y$ whose curvature
is a positive current \cite{d2}.%
\footnote{the somewhat confusing terminology of positive currents actually means
that $dd^{c}\phi$ is allowed to be a \emph{semi-}positive form on
the set where $\phi$ is smooth.%
} The domain $X$ in $Y$ will be called a \emph{polarized domain}
if the metric $\phi$ on $L$ is smooth on $X$ and it will be called
a \emph{polarized domain with compatible curvatures} if any of the
assumptions in section \ref{sub:Assumptions-on-the} below are satisfied.
\footnote{Since a polarization usually refers to a \emph{positive} line bundle
$L,$ the term semi-polarized would perhaps be more appropriate.%
} 

Fixing an Hermitian metric two-form $\omega$ on $X$ (with associated
volume form $\omega_{n})$ the Hilbert space $\mathcal{H}(Y,L^{k})_{X}$
is defined as the space $H^{0}(Y,L^{k})$ with the norm obtained by
restriction of the global norm on \cite{gri} to $X:$

\begin{equation}
\left\Vert \alpha_{k}\right\Vert _{k\phi}^{2}:=\left\Vert \alpha_{k}\right\Vert _{k\phi X}^{2}(=\int_{X}\left|f_{k}\right|^{2}e^{-k\phi(z)}\omega_{n}),\label{eq:norm restr}\end{equation}
using a suggestive notation in the last equality (compare formula
\ref{eq:ptwise norm}). If $\eta$ is a form we will write $\eta_{p}:=\eta^{p}/p!,$
so that the volume form on $X$ may be written as $\omega_{n}.$ The
induced volume form on $\partial X$ will be denoted by $d\sigma.$
If $Z$ is a submanifold, then $[Z]$ will denote the corresponding
current, i.e. $([Z],\eta):=\int_{Z}\eta$ for any test form $\eta.$
Moreover, given a real $(1,1)-$form $\eta$ on $Y$ we will denote
by $\{\eta\}_{y}$ the corresponding (using the metric form $\omega)$
Hermitian linear operator on $T^{1,0}(Y)_{y}$ (or on some specified
subbundle).

\subsection{\label{sub:Assumptions-on-the}Assumptions for {}``compatible curvatures''}

At least one of the following three assumptions are assumed to be
satisfied for $X$ to be a \emph{polarized domain with compatible
curvatures} (compare \cite{berm4}). The assumptions all have in common
the condition that the slope function (see formula \ref{eq:def of mu})
is constant:\begin{equation}
T\equiv C\,\textrm{on\,}\partial X,\label{eq:ass on T general}\end{equation}
for some choice of the defining function $\rho$ of $\partial X.$

\subsubsection*{Assumption 1}

The defining function $-\rho$ of the pseudoconvex manifold $Y-X$
may be chosen to be smooth with $d\rho\neq0$ and $dd^{c}(-\rho)>0$
in $(Y-X)-Z,$ where $Z$ is either a point or an irreducible divisor
in $Y-X.$%
\footnote{i.e. a (possibly singular) connected compact closed complex submanifold
of codimension one in $Y-X.$ Then the integration current $[Z]$
is well-defined. \cite{gri} %
} Moreover, on any regular sublevelset of $\rho$ the slope function
$T$ in \ref{eq:def of mu} (defined by replacing $\partial X$ with
the sublevelset of $\rho)$ is constant, i.e. \begin{equation}
T\textrm{\, is\, a\, function\, of\,\,$\rho.$}\label{eq:ass T is f of rho}\end{equation}
 If $Z$ is a point it is assumed that $dd^{c}(-\rho)>0$ on all of
$Y-X.$ If $Z$ is an irreducible divisor it is assumed that $T$
is bounded from above on $Y-X,$ that \begin{equation}
\int_{Z}c_{1}(L)^{n-1}=0.\label{eq:ass on pairing}\end{equation}
and that \begin{equation}
dd^{c}(-\rho)=[Z]+\beta,\label{eq:ass on sing}\end{equation}
in the sense of currents on $Y-X,$ where $\beta$ is a semi-positive
smooth form.

\subsubsection*{Assumption 2}

Suppose that $n\geq2$ (the dimension of $X)$ and that $L$ is holomorphically
trivial on $Y-X.$ Then the fiber metric $\phi$ on $L$ may be identified
with a function on $Y-X$ and it as assumed that \[
\phi=-\rho\]
on $Y-X.$ In this case the form $\mu$ in formula \ref{eq:def of mu}
is simply given by \[
\mu=(dd^{c}\phi)_{n-1}\wedge d^{c}\phi/n.\]

\subsubsection*{Assumption 3}

Suppose that $n\geq3,$ that $Y-X$ is a Stein manifold and that \begin{equation}
dd^{c}\phi=-fdd^{c}\rho\label{eq:ass conf equi}\end{equation}
along the holomorphic tangentbundle of $\partial X$ for some non-negative
function $f$ on $\partial X$.

\subsection{\label{sub:General-properties-of}General properties of Bergman kernels}

Let $(\psi_{i})$ be an orthonormal base for a given Hilbert space
structure on the space $H^{0}(Y,L^{k}),$ which in this paper always
will be the Hilbert space $\mathcal{H}_{k}(Y,L^{k})_{X}$. The \emph{Bergman
kernel} of the Hilbert space $H^{0}(Y,L)$ is defined by \[
K^{k}(x,y)=\sum_{i}\psi_{i}(x)\otimes\overline{\psi_{i}(y)}.\]
 Hence, $K^{k}(x,y)$ is a section of the pulled back line bundle
$L^{k}\boxtimes\overline{L}^{k}$ over $Y\times Y.$ For a fixed point
$y$ we identify $K_{y}^{k}(x):=K^{k}(x,y)$ with a section of the
hermitian line bundle $L^{k}\otimes L_{y}^{k},$ where $L_{y}$ denotes
the line bundle over $Y,$ whose constant fiber is the fiber of $L$
over $y,$ with the induced metric. The definition of $K^{k}$ is
made so that $K^{k}$ satisfies the following reproducing property
\begin{equation}
\alpha(y)=(\alpha,K_{y}^{k})_{k\phi}\label{(I)repr property}\end{equation}
\footnote{We are abusing notation here: the scalar product $(\cdot,\cdot)_{k\phi}$
on $H^{0}(Y,L^{k})$ determines a pairing of $K_{y}^{k}$ with any
element of $H^{0}(Y,L^{k}),$ yielding an element of $L_{y}^{k}.$ %
}for any element $\alpha$ of $H^{0}(Y,L^{k}),$ which also shows that
$K^{k}$ is well-defined. In other words $K^{k}$ is the integral
kernel of the orthogonal projection onto $H^{0}(Y^{k},L)$ in $L^{2}(Y,L^{k}).$
The restriction of $K^{k}$ to the diagonal is a section of $L^{k}\otimes\overline{L}^{k}$
and we let $B^{k}(x)=\left|K^{k}(x,x)\right|_{k\phi}(=\left|K^{k}(x,x)\right|e^{-k\phi(x)})$
be its point wise norm: \[
B^{k}(x)=\sum_{i}\left|\psi_{i}(x)\right|_{k\phi}^{2}.\]
We will refer to $B^{k}(x)$ and $B^{k}1_{X}\omega_{n}$ as the \emph{Bergman
function} and \emph{Bergman measure of} $H^{0}(Y,L^{k}).$ Note that
the Bergman measure only depends on the {}``restriction'' of the
metric $\phi$ to the domain $X.$ The following extremal property
holds:\begin{equation}
B^{k}(x)=\sup\left|\alpha_{k}(x)\right|_{k\phi}^{2},\label{(I)extremal prop of B}\end{equation}
 where the supremum is taken over all $L^{2}-$normalized elements
$\alpha_{k}$ of $H^{0}(Y,L^{k}).$ An element realizing the extremum,
is called an \emph{extremal at the point $x$} and is determined up
to a complex constant of unit norm. Given such an extremal $\alpha$
the following basic relation holds \cite{berm2}: \begin{equation}
\left|K^{k}(x,y)\right|_{k\phi}^{2}=\left|\alpha_{k}(y)\right|_{k\phi}^{2}B^{k}(x)\label{eq:K as extremal}\end{equation}

\section{Examples\label{sec:Model-examples} and {}``counter examples''}

In this section we will consider various classes of polarized pseudoconcave
domains. Some {}`` counter examples'' will also be presented, showing
that the main results in this article may not hold if the assumptions
in section \ref{sub:Assumptions-on-the} are relaxed.

\subsection{Domains in projective space and polynomials in $\C^{n}$}

Let $Y$ be the $n-$dimensional projective space $\P^{n}$ and let
$L$ be the hyperplane line bundle $\mathcal{O}(1)$. Then $H^{0}(Y,L^{k})$
is the space of homogeneous polynomials in $n+1$ homogeneous coordinates
$Z_{0},Z_{1},..Z_{n}$ \cite{gri}. The Fubini-Study metric $\phi_{FS}$
on $\mathcal{O}(1)$ may be suggestively written as $\phi_{FS}(Z)=\ln(\left|Z\right|^{2})$
and the Fubini-Study metric $\omega_{FS}$ on $\P^{n}$ is the normalized
curvature form $dd^{c}\phi_{FS},$ which is hence invariant under
the standard action of $SU(n+1)$ on $\P^{n}.$ We may identify $\C^{n}$
with the {}``affine piece'' $\P^{n}-H_{\infty}$ where $H_{\infty}$
is the {}``hyperplane at infinity'' in $\C^{n}$ (defined as the
set where $Z_{0}=0).$ In terms of the standard trivialization of
$\mathcal{O}(1)$ over $\C^{n},$ the space $H^{0}(Y,L^{k})$ may
be identified with the space of polynomials $p_{k}(\zeta)$ in $\C_{\zeta}^{n}$
of total degree at most $k$ and the fiber metric $\phi$ may be identified
with a (plurisubharmonic) function in $\C^{n}.$ 

The basic example of a pseudoconcave domain $X$ is obtained as the
complement in $\P^{n}$ of the unit-ball in $\C^{n}$ (i.e. we may
take $\rho=-\left|\zeta\right|^{2}+1$ in a neighbourhood of $\partial X).$
The norm \ref{eq:norm restr} on the Hilbert space $\mathcal{H}_{k}(Y,L^{k})_{X}$
may in this case be expressed as

\[
\left\Vert p_{k}\right\Vert _{k\phi}^{2}:=\int_{\left|z\right|>1}\left|p_{k}(\zeta)\right|^{2}e^{-k\phi(\zeta)}(\omega_{FS})_{n}\]

\begin{example}
The canonical Fubini-Study metric on $\mathcal{O}(1)$ corresponds
to the choice $\phi(\zeta)=\ln(1+\left|\zeta\right|^{2}).$ In this
case the (normalized) curvature of $L$ is the Fubini-Study metric
on $\P^{n}$ and the form $\mu$ in formula \ref{eq:def of mu} is
a multiple of the standard volume form on the $2n-1-$sphere.
\end{example}
Next, we will consider the basic example when $\phi$ has a singularity
in $Y-X.$ 

\begin{example}
\label{ex: lnplus}Let $\phi(\zeta)=\ln(\left|\zeta\right|^{2}).$
Then $\phi$ is smooth outside the origin. In particular it is smooth
on $X.$ Note that the $n$th exterior power of the curvature of $\phi$
vanishes outside the origin (i.e. the complex Monge-Ampere of $\phi$
vanishes there). Again, the form $\mu$ is a multiple of the standard
volume form on the $2n-1-$sphere.
\end{example}
In the following section generalizations of the case when $X$ is
the complement of the unit-ball are considered (compare remark \ref{rem:blow up as bundle}).

\subsection{Disc bundles\label{sub:Discbundles}}

Let $Z_{+}$ be a closed compact complex manifold of dimension $n-1$
and let $(F,\phi_{F})$ and $(G,\phi_{G})$ be Hermitian holomorphic
line bundle over $Z_{+}$ with positive curvature. Then $X$ is defined
as the pseudoconcave domain obtained as the unit discbundle in the
total space of $F$ and $Y$ as the $\P{}^{1}-$bundle over $Z$ obtained
by fiberwise {}``adding the point at infinity'' to $F,$ i.e. by
adding a divisor $Z_{-}$ at infinity. %
\footnote{i.e. $Y$ is the fiber-wise projectivization of the bundle $F\oplus\underline{\C}$,
where $\underline{\C}$ is the trivial line bundle over $Z_{+}.$
The coordinate along $\underline{\C}$ determines a section of $\mathcal{O}_{\P(F\oplus\underline{\C)}}(1)$
whose zero-set is $Z_{-}.$ %
}Hence, locally \[
X=\{ h=\left|w\right|^{2}\exp(-\phi_{F}(z))\leq1\}\]
 (where $z$ is a coordinate along $Z$ and $w$ is a coordinate along
the fibers of $F).$ We will assume that the {}``slope between $dd^{c}\phi_{G}$
and $dd^{c}\phi_{F}$'':\begin{equation}
S=\sup\left\{ t\geq0:\,(dd^{c}\phi_{G}-tdd^{c}\phi_{F})_{z}\geq0\right\} \equiv S_{0}\in\Z,\label{eq:ass S is indep}\end{equation}
i.e. that $S$ above is independent of the point $z$ in $Z$ and
that $S\in\Z.$ The line bundle $L$ over $Y$ is now defined as \[
L:=\pi_{F}(G)\otimes[Z_{-}]^{S_{0}}\]
where $[Z_{-}]$ now denotes the \emph{line bundle} over $Y$ corresponding
to the divisor $Z_{-}$(do that $c_{1}([Z_{-}])$ is represented by
the \emph{current} $[Z_{-}],$ using the notation introduced in section
\ref{sec:Setup}) Hence, $L$ is isomorphic to $\pi_{F}^{*}(G)$ over
the domain $X.$ As will be seen below $L$ satisfies the assumption
\ref{eq:ass on pairing} (with $Z=Z_{-}).$ Moreover, the defining
function $\rho:=\ln h$ for the boundary of $X$ satisfies the assumption
\ref{eq:ass on sing} (with $\beta=dd^{c}\phi_{F}$). Let us now consider
two different metrics on $L:$

\begin{example}
Let $\phi_{L}(z,w):=\pi_{G}^{*}\phi_{G}+\ln(1+e^{S\rho})$ on $Y-Z_{-}$
(smoothly extended as a metric on $L$ over a neighbourhood of $Z_{-}).$
Then $\phi_{L}$ has positive curvature on $X-Z_{-}$ and $(dd^{c}\phi_{G})^{n-1}=0$
on $Z_{-}$ precisely when assumption \ref{eq:ass S is indep} holds.
Indeed, a direct calculation (compare the proof of lemma \ref{lem:mod met})
gives \[
dd^{c}\phi_{L}=\pi_{F}^{*}(dd^{c}\phi_{G}-s(\rho)dd^{c}\phi_{F})+dt\wedge\pi_{F}^{*}d^{c}\rho,\]
 where $s(\rho)=\frac{\partial}{\partial\rho}\ln(1+e^{S\rho})$ is
strictly increasing, mapping $[-\infty,\infty]$ to $[0,S].$ Note
that the slope function $T$ (see formula \ref{eq:def of mu}) becomes
a function of $\rho:$ $T(\rho)=S-s(\rho)$ in this case. Hence, all
assumptions in 1 in section \ref{sub:Assumptions-on-the} are satisfied
and $(X,\phi_{L})$ is thus a polarized domain with compatible curvatures
in $(Y,L).$
\end{example}
The following example may be obtained as a limit of variants of the
previous one:

\begin{example}
\label{exa:pull-back bd}Let $\phi_{L}(z,w):=\pi_{F}^{*}\phi_{G}(z)$
on $X.$ Then $\phi_{L}$ extends to a singular metric on $L$ over
$Y$ with positive curvature (in the sense of currents) by setting
$\phi_{L}(z,w):=\phi_{G}(z)+\rho$ on $Y-X.$ Indeed, $\phi_{L}$
may be obtained as the limit \[
\phi_{L}(z,w):=\lim_{k\rightarrow\infty}(\pi_{G}^{*}\phi_{G}+k^{-1}\ln(1+e^{Sk\rho}))\]
Note that in this example $S=T.$
\end{example}
\begin{rem}
\label{rem:blow up as bundle}Setting $Z_{+}=\P^{n-1}$ and $F=G=\mathcal{O}(1)$
gives the case considered in the previous section, i.e. when $X$
is the complement of the unit-ball in $\P^{n}.$ Indeed, the base
$Z_{+}$ of the fibration above may be identified with the hyperplane
at infinity in $\C^{n}$ and the fibers correspond to lines through
the origin in $\C^{n}$. Moreover, $Y$ corresponds to the blow-up
$\widetilde{\P^{n}}$ at the origin of $\P^{n}.$ A local isomorphism
between $Y$ and $\widetilde{\P^{n}}$ is obtained by setting $w=\zeta_{n}^{-1}$
and $z=(\zeta_{1}/\zeta_{n},...,\zeta_{n-1}/\zeta_{n}).$ Also note
that $Z_{-}$ corresponds to the exceptional divisor $E$ over the
origin and $L\approx\pi^{*}\mathcal{O}_{\P^{n}}(1),$ where $\pi$
is the blow-down map from $\widetilde{\P^{n}}$ to $\P^{n}.$ Hence,
the space $H^{0}(Y,L^{k})$ may be identified with $H^{0}(\P^{n},\mathcal{O}_{\P^{n}}(1)^{k})$
in this particular case.
\end{rem}

\subsection{{}``\label{sub:Counter-examples}Counter examples''}

\subsubsection{Constant vs. non-constant slope $T$}

Now assume that the base $Z_{+}$ is a product of complex curves:
$Z_{+}=Z_{1}\times Z_{2},$ $\dim Z_{i}=1.$ All other objects are
assumed to decompose accordingly: $F=F_{1}\otimes F_{2}$ etc. Then
the curvature form $dd^{c}\phi_{F}$ may be identified with a pair
of functions (the eigenvalues of $dd^{c}\phi_{F}):$ $(\lambda_{F_{1}}(z_{1}),\lambda_{F_{2}}(z_{2}))$
and similarly for $dd^{c}\phi_{G}.$ For concreteness, consider the
case when \[
dd^{c}\phi_{G}\leftrightarrow(2,2),\,\,\, dd^{c}\phi_{F}\leftrightarrow(1+\epsilon_{1}(z_{1}),2+\epsilon_{2}(z_{2})),\]
 where $\left|\epsilon_{1}(z_{1})\right|<1,$ $\left|\epsilon_{2}(z_{2})\right|<2$
and the integral of $\epsilon_{i}(z_{i})$ over $Z_{i}$ vanishes.
Hence, $(\lambda_{F_{1}}(z_{1}),\lambda_{F_{2}}(z_{2}))$ corresponds
to a deformation (with positive curvature) of a metric on $F$ with
constant curvature. Now, if $\epsilon_{1}\equiv0,$ then the slope
function $T(z)$ is clearly \emph{constant} ($=1),$ as long as $\left|\epsilon_{2}(z_{2})\right|\leq1.$
This means that the corresponding disc bundle $X$ with the hermitian
line bundle $L$ defined as in example \ref{exa:pull-back bd} is
hence a polarized domain with compatible curvatures (since it satisfies
the assumption 1 in section \ref{sub:Assumptions-on-the}). But if
$\epsilon_{1}\neq0$ somewhere and if $\epsilon_{2}(z_{2})\leq\epsilon_{1}(z_{1})+1,$
then a short calculation gives\[
T(z)=2/(1+\epsilon_{1}(z_{1}))\]
 and $T$ is hence \emph{non-constant.} In fact, it can be shown that
\emph{the main results of this paper hold if and only if $T(z)$ is
constant in these examples} (and certain more general examples) \emph{\cite{berm5},}
which makes the assumption \ref{eq:ass on T general} quite natural. 

Finally, consider the case when $Z_{+}$ is a complex curve, so that
$X$ is a domain in a complex surface (i.e. $n=2).$ Then $T(z)=dd^{c}\phi_{G}(z)/dd^{c}\phi_{F}(z).$
Now suppose that $Z_{+}=\P^{1}$ and $F=G=\mathcal{O}(1)$ with $\phi_{G}(z)=\ln(1+\left|z\right|^{2})$
and $\phi_{F}(z)=\ln(a+\left|z\right|^{2})$ (compare remark \ref{rem:blow up as bundle}).
Then $X$ corresponds to the exterior of the {}``ellipse'' $\{(\zeta_{1},\zeta_{2}):\,\left|\zeta_{1}\right|^{2}+a\left|\zeta_{2}\right|^{2}=1\}$
in $\C^{2}.$ Hence, $T(z)$ is constant precisely when $a=1,$ i.e.
when $X$ is the exterior of the unit-ball in $\C^{2}.$ This example
also shows the need to assume that $n>2$ in assumption 3 in section
\ref{sub:Assumptions-on-the}. Indeed, when $n=2$ the assumption
\ref{eq:ass conf equi} always holds. Similarly, the example also
 shows the need to assume extension properties of $\rho$ and $(L,\phi),$
as in the assumptions 1 and 2 (at least when $n=2),$ since \ref{eq:ass on T general}
always holds when $n=2$ (for a suitable choice of $\rho).$

\subsubsection{Vanishing vs. non-vanishing of $c(L)^{n-1}\cdot[Z]$}

The following example illustrates the need for the condition \ref{eq:ass on pairing}
in the assumption 1 in section \ref{sub:Assumptions-on-the}. Consider
the situation in remark \ref{rem:blow up as bundle}, but replace
$L$ with the line bundle $L=\pi^{*}\mathcal{O}_{\P^{n}}(2)\otimes[E]^{-1},$
Then $L$ is a positive line bundle and satisfies all the assumptions
in 1 in section \ref{sub:Assumptions-on-the}, except \ref{eq:ass on pairing}.
Indeed, $L$ corresponds to $\pi_{F}(G)\otimes[Z_{-}]^{S_{0}-1}$
in section (with $G=\mathcal{O}_{\P^{n-1}}(2)$ and $F=\mathcal{O}_{\P^{n-1}}(1))$.
However, the assumption \ref{eq:ass on pairing} fails, since\[
c(L)^{n-1}\cdot[E]=0+0+...+0+[E]^{n}\neq0\]
Note that $L$ has a natural singular metric with curvature form $\pi^{*}c(\mathcal{O}(2))+E,$
where $c(\mathcal{O}(2))$ denotes the curvature form of a fixed smooth
hermitian metric on $\mathcal{O}(2)$ with positive curvature, so
that $L$ is isomorphic to $\mathcal{O}(2)$ over $X$ (as Hermitian
holomorphic line bundles). Moreover, there is a strict inclusion\begin{equation}
H^{0}(\widetilde{\P^{n}},L^{k})\hookrightarrow H^{0}(\widetilde{\P^{n}},\pi^{*}\mathcal{O}_{\P^{2}}(2)^{k}),\label{eq:inclus}\end{equation}
 where the image is the subspace in $H^{0}(\widetilde{\P^{n}},\pi^{*}\mathcal{O}_{\P^{2}}(2)^{k})$
of all sections vanishing along $E$ to order $k$ (i.e. the image
may be identified with the subspace of all polynomials in $\C^{n}$
of total degree $m$, where $k<m\leq2k).$ It follows that the main
results of this paper do not apply to $(\widetilde{\P^{n}},L),$ since
they imply that formula \ref{pro:Morse eq} holds, where the left
hand side in the formula only depends on the restriction of the curvature
of $L$ to $X.$ Indeed, the formula hold does hold for $(\widetilde{\P^{n}},\pi^{*}\mathcal{O}_{\P^{2}}(2)),$
since it satisfies all the assumptions in 1. Hence, by \ref{eq:inclus}
it cannot hold for $(\widetilde{\P^{n}},L).$

\section{\label{sec:Morse-(in)equalities-and}Morse (in)equalities and model
Bergman kernels}

In this section we will mainly recall some point-wise estimates for
the Bergman functions of the space $\mathcal{H}_{k}(Y,L^{k})_{X}$
obtained in \cite{berm1,berm3}. Such point-wise estimates were referred
to as \emph{local holomorphic Morse inequalities} in \cite{berm1}
in the more general context of harmonic $(0,q)-$form with values
in $L^{k}$.%
\footnote{The case of holomorphic sections is considerably more elementary than
the general case. The main difference is that there is no need for
a special sequence of metrics on $X$ as in \cite{berm3} and that
subelliptic estimates may be replaced by the submean property of holomorphic
functions.%
} After integration the latter estimates yield bounds for the asymptotic
growth of the Dolbeault cohomology groups with values in $L^{k}.$
The latter bounds were first obtained by Demailly \cite{d1} in the
context of closed manifolds, who called them \emph{holomorphic} \emph{Morse
inequalities} in analogy with the classical case of Morse inequalities
for the De Rham cohomology groups of a real manifold (compare Witten's
approach in \cite{wi}).

\subsection{\label{sub:int reg}Morse inequalities in the interior region}

For a fixed $k$ the interior region is defined by the inequality
$\rho\leq-1/\ln k.$%
\footnote{in the following $\ln k$ could be replaced with any sequence $R_{k}$
tending to infinity at the order $O(k^{\epsilon})$ where $\epsilon$
is a sufficiently small positive number. Note that in \cite{berm3}
$R_{k}=R$ and the limit when first $k$ and then $R$ tend to infinity
was considered. In this paper a slightly more precise control in the
{}``boundary region'' (see the appendix) will allow us to let $R$
depend on $k,$ hence simplifying the notation a bit.%
} In \cite{berm3} it was shown that the Bergman function $B^{k}(x)$
may be estimated in terms of \emph{model Bergman functions.} The model
Bergman function $B^{0}$ associated to an interior point $x$ is
obtained by replacing the manifold $X$ with $\C^{n}$ and the line
bundle $L$ with the constant curvature line bundle over $\C^{n}$
obtained by freezing the curvature of $L$ at the point $x.$ Since
$\C^{n}$ is non-compact all sections are assumed to have finite $L^{2}-$norm.
More concretely, one may always arrange so that locally around the
fixed point $x,$ \begin{equation}
\phi(z)=\sum_{i=1}^{n}\lambda_{i}\left|z_{i}\right|^{2}+...,\,\,\,\omega(z)=\frac{i}{2}\sum_{i=1}^{n}dz_{i}\wedge\overline{dz_{i}}+....\label{eq: Intro metrics}\end{equation}
 where the dots indicate lower order terms and the leading terms are
called \emph{model metrics.} Hence, the corresponding model $L^{2}-$norm
on $\C^{n}$ is given by\begin{equation}
\int_{\C^{n}}\left|\alpha(z)\right|^{2}e^{-\sum_{i=1}^{n}\lambda_{i}\left|z_{i}\right|^{2}},\label{eq: Intro model norm}\end{equation}
integrating with respect to the Euclidean measure on $\C^{n}.$ 

Denote by $F_{k}$ the holomorphic scaling map \[
F_{k}(z)=(z/k^{1/2})\]
and let $\alpha^{(k)}(z):=(F_{k}^{*}\alpha_{k})(z).$ By the proof
of theorem 1.1. in \cite{berm1} (see also \cite{berm2} for a simple
argument based on the submean property of holomorphic functions) the
following point-wise bound holds in the interior region: \begin{equation}
\limsup_{k}k^{-n}\left|\alpha^{(k)}(z)\right|_{k\phi}^{2}/\left\Vert \alpha_{k}\right\Vert _{k\phi F_{k}(D_{\ln k})}^{2}\leq B^{0}(0),\label{eq:local Morse S int}\end{equation}
where $D_{lnk}$ denotes a polydisc of radius $\textrm{ln\,}k.$ In
particular, by the extremal property \ref{(I)extremal prop of B}
of $B^{k}(x):$ 

\begin{equation}
\limsup_{k}k^{-n}B^{(k)}(z)\leq B^{0}(0)\label{eq:Intro local weak pointwise}\end{equation}
 Moreover, the model Bergman function is explicitly given by\begin{equation}
B^{0}\omega_{n}=(dd^{c}\phi)_{n}\label{eq:model B int}\end{equation}
 Hence, the full model Bergman kernel is given by \begin{equation}
K^{0}=\det(dd^{c}\phi_{0})e^{\phi_{0}},\label{eq:model K int}\end{equation}
using the suggestive notation $\phi_{0}=\phi_{0}(z,z')=\sum_{i=1}^{n}\lambda_{i}\overline{z_{i}}z_{i}'.$

\subsection{\label{sub:mid reg}Morse inequalities in the middle region}

The middle region is defined by the inequalities $-1/\ln k\leq\rho\leq-\ln k/k.$
As shown in \cite{berm3} (section 5.2)\begin{equation}
\lim_{k}\int_{-1/\ln k\leq\rho\leq-\ln k/k.}k^{-n}B^{k}\omega_{n}=0\label{eq:ineq int middle}\end{equation}

\subsection{\label{sub:bd reg}Morse inequalities in the boundary region}

The boundary region is defined by the inequalities $-\ln k/k\leq\rho\leq0$
and is diffeomorphic to the product $\partial X\times[-\ln k/k,0].$
Fix a point $\sigma$ in $\partial X$ and take local holomorphic
coordinates $(z,w),$ where $z$ is in $\C^{n-1}$ and $w=u+iv.$
By an appropriate choice we may assume that the coordinates are orthonormal
at $0$ and that \begin{equation}
\rho(z,w)=v+\sum_{i=1}^{n-1}\mu_{i}\left|z_{i}\right|^{2}+O(\left|(z,w)\right|^{3})=:\rho_{0}(z,w)+O(\left|(z,w)\right|^{3}).\label{II local expr for def funct}\end{equation}
 In a suitable local holomorphic trivialization of $L$ close to the
boundary point $\sigma,$ the fiber metric may be written as \[
\phi(z,w)=\sum_{i,j=1}^{n-1}\lambda_{ij}z_{i}\overline{z_{j}}+O(\left|w\right|)O(\left|z\right|)+O(\left|w\right|^{2})+O(\left|(z,w)\right|^{3}),\]
 where the leading terms are called the model fiber metric and denoted
by $\phi_{0}.$ The model Bergman function $B^{0}$ and kernel $K^{0}$
associated to the fixed point $\sigma$ are the ones obtained from
the Hilbert space of all holomorphic functions $\alpha$ on the model
domain $X_{0}$ (with defining function $\rho_{0})$ which are square
integrable with respect to the model norm \[
\int_{X_{0}}\left|\alpha(z,w)\right|^{2}e^{-\phi_{0}(z)},\]
integrating with respect to the Euclidean measure on $\C_{z,w}^{n}.$ 

Denote by $F_{k}$ the holomorphic scaling map \begin{equation}
F_{k}(z,w)=(z/k^{1/2},w/k),\label{eq:scal map bd}\end{equation}
 so that \[
X_{k}=F_{k}(D_{\ln k})\bigcap X\]
 is a sequence of decreasing neighborhoods of the boundary point $\sigma,$
where $D_{\ln k}$ denotes the polydisc of radius $\ln k$ in $\C^{n}.$
Note that \begin{equation}
F_{k}^{-1}(X_{k})\rightarrow X_{0},\label{eq:conv of domains}\end{equation}
 in a suitable sense. On $F_{k}^{-1}(X_{k})$ we have the \emph{scaled}
fiber metric $F_{k}^{*}k\phi$ that tends to the model metric $\phi_{0}$
on the model domain $X_{0},$ when $k$ tends to infinity. 

It follows from the proof of proposition 5.5 in \cite{berm3} that
\begin{equation}
\limsup_{k}k^{-(n+1)}\left|\alpha^{(k)}(z,w)\right|_{k\phi}^{2}/\left\Vert \alpha_{k}\right\Vert _{k\phi X_{k}}^{2}\leq B^{0}(z,w),\label{eq:Morse ineq S bd}\end{equation}
where the estimate is uniform on $F_{k}^{-1}(X_{k}).$ Moreover, the
left hand side above is uniformly bounded by a constant on any polydisc
of fixed radius in $\C^{n}$ (even without intersecting with $X)$
and on $X_{0}$ dominated by an $L^{1}-$function on {}``half-rays''
(see the appendix). In particular, \begin{equation}
\begin{array}{lrcl}
(i) & \limsup_{k}k^{-(n+1)}B^{(k)}(z,w) & \leq & B^{0}(z,w)\\
\end{array}\label{eq:ineq B bd}\end{equation}
The model Bergman function at the fixed point $\sigma$ in $\partial X$
is given by\begin{equation}
B_{X_{0}}(z,w)\omega_{n}=\int_{0}^{T}e^{t\rho_{0}(z,w)}t(dd^{c}\phi_{0}+tdd^{c}\rho_{0})_{n-1}\wedge d^{c}\rho_{0})dt,\label{eq:model B bd}\end{equation}
where $T$ is the slope function in formula \ref{eq:def of mu}. In
particular, performing the fiber-integral (i.e the push forward) over
$\rho$ gives\begin{equation}
\int_{\rho=-\infty}^{0}B^{0}\omega_{n}=\mu,\label{eq:model B bd fiberint}\end{equation}
where $\mu$ as the $2n-1$ form defined by formula \ref{eq:def of mu}.
Similarly, the corresponding model Bergman kernel is given by the
following integral formula\begin{equation}
K^{0}=\frac{1}{4\pi}\frac{1}{\pi}\int_{0}^{T}e^{t\rho_{0}+\phi_{0}}t\textrm{det}(dd^{c}\phi_{0}+tdd^{c}\rho_{0})dt,\label{eq:model K bd}\end{equation}
 using a suggestive notation as in formula \ref{eq:model K int} above.
Equivalently, we have the suggestive formula\[
K^{0}=\frac{1}{4\pi}\frac{1}{\pi}\textrm{det}(dd^{c}\rho_{0})e^{\phi_{0}}P(\frac{\partial}{\partial\rho_{0}})\frac{\partial}{\partial\rho_{0}}(\frac{e^{T\rho_{0}}-1}{\rho_{0}}),\]
 where $P$ is the characteristic polynomial of the operator $\{ dd^{c}\phi\}\{-dd^{c}\rho\}^{-1}$
where the operators act on $T^{1,0}(\partial X)_{x}$and $(P(\frac{\partial}{\partial\rho_{0}})$
denotes the corresponding differential operator with constant coefficients).
Note that $T$ is the minimal eigenvalue of $\{ dd^{c}\phi\}\{-dd^{c}\rho\}^{-1}.$

\subsection{Morse equalities}

Integrating the Bergman function $B^{k}$ over $X$ gives, using the
point-wise bounds in the different regions the following Morse inequalities
for any line bundle $L$ over $X$ (compare \cite{berm3}) \[
\mathcal{\dim H}_{k}(Y)_{X}\leq k^{n}(\int_{X}(dd^{c}\phi)_{n}+\int_{\partial X}\mu)+o(k^{n}),\]
where the form $\mu$ is given by formula \ref{eq:def of mu}.

In the case when $X$ is a polarized pseudoconcave domain (section
\ref{sub:Notation-and-setup}) the previous inequality becomes an
equality:

\begin{prop}
\label{pro:Morse eq}Consider the Hilbert space $\mathcal{H}_{k}(Y,L^{k})_{X}$
associated to the polarized pseudoconcave domain $X.$ Then \begin{equation}
\lim_{k}k^{-n}\mathcal{\dim H}_{k}(Y)_{X}=(\int_{X}(dd^{c}\phi)_{n}+\int_{\partial X}\mu)\label{eq:Morse eq}\end{equation}

\end{prop}
\begin{proof}
\emph{The assumption 1 holds:} If the extension of the fixed fiber
metric $\phi$ to $Y$ (that will also be denoted by $\phi$ in the
sequel) is smooth, then the equality in \ref{eq:Morse eq} was obtained
in \cite{berm3} (section 7.1), but stated incorrectly there without
any assumptions on the slope function $T$ in formula \ref{eq:def of mu}).
We will next repeat the argument in \ref{eq:Morse eq} and point out
the corrections appearing in \cite{berm4}. First, since $L$ is a
semi-positive line bundle the left hand side above is given by, using
for example Demailly's strong Morse inequalities \cite{d1}: \[
\int_{Y}c^{1}(L)^{n}=\int_{X}c^{1}(L)^{n}+\int_{Y-X}c^{1}(L)^{n}.\]
Next, Stokes theorem is used to show that the integral over $Y-X$
coincides with the boundary integral in \ref{eq:Morse eq} (also using
that $dd^{c}\phi$ is a smooth semi-positive form representing $c^{1}(L)$)
(proposition 0.1 in \cite{berm4}). Note that we may even assume that
$dd^{c}\phi$ is strictly positive by adding $\omega/m$ to $dd^{c}\phi,$
where $\omega$ is the curvature form of any fixed positive line bundle
$A$ on $X$ and in the end letting the positive number $m$ tend
to zero in the integrals. 

The case when the extension $\phi$ is singular may be reduced to
the smooth case. Indeed, for a fixed $m$ the positivity of $L^{m}\otimes A$
allows us to find a sequence $\Phi_{j}$ of smooth metrics on $L^{m}\otimes A$
with semi-positive curvature decreasing to $\Phi,$ where $dd^{c}\Phi=mdd^{c}\phi+\omega$
(using the regularization results of Demailly - see the appendix in
\cite{g-z} for a direct proof). This gives the Morse equality \ref{eq:Morse eq}
with $\phi$ replaced by $\Phi_{j}$. Finally, letting first $j$
and then $m$ tend to infinity then proves the proposition in the
singular case, since the operator that maps $\phi$ to $(dd^{c}\phi)^{p}$
is continuous when applied to a locally bounded decreasing sequence
of smooth plurisubharmonic functions, as first shown by Bedford-Taylor
(see \cite{kl}).

\emph{The assumption 2 holds:} in this case proposition 0.1 in \cite{berm4}
may be replaced by a direct application of Stokes theorem (compare
\cite{berm3}, section 7.1).

\emph{The assumption 3 holds:} under the assumption \ref{eq:ass conf equi}
the strong Morse inequalities obtained in \cite{berm3} give \ref{eq:Morse eq}
with $\mathcal{\dim H}_{k}(Y)_{X}(=H^{0}(Y,L^{k}))$ replaced by $H^{0}(X,L^{k})$.
But if $Y-X$ is a Stein manifold a standard extension argument then
shows that $H^{0}(X,L^{k})=H^{0}(Y,L^{k}),$ using that the Dolbeault
cohomology group $H_{cpt}^{0,1}(M,F)$ for compactly supported forms
is trivial for any line bundle $F$ on a Stein manifold $M$ of dimension
at least three (setting $M=Y-X$ and $F=L^{k}).$\cite{d2}
\end{proof}
The following simple generalization of lemma 2.2 in \cite{berm2}
will be used to convert the Morse equalities from the previous proposition
to equalities for the \emph{scaled} Bergman functions and kernels.

\begin{lem}
\label{lem:integrat theory}Assume that $(M,\nu)$ is a manifold with
a smooth volume form $\nu$ and $V$ an open subset of $M.$ For each
fixed parameter $u$ defined in a bounded set $U$ of Euclidean $\R^{N}$
let $\Phi_{k,u}$ be a diffeomorphism from $V$ onto $\Phi_{k,u}(V)$
such that the Jacobian of $\Phi_{k,u}$ converges uniformly to $1$,
when $k$ tends to infinity. Let $f$ and $f_{k}$ be functions in
$L^{1}(V,\nu)$ with compact support in $V$ and such that $\textrm{supp}f_{k}\subset\Phi_{k,u}(V)$
and such that the sequence $f_{k}$ is dominated by a function in
$L^{1}(V,\nu)$. Moreover, assume that\[
(\textrm{i})\lim_{k}\int_{V}f_{k}d\nu=\int_{V}fd\nu\,\,\,\,\textrm{and}\,\,\,\,(\textrm{ii})\,\limsup f_{k}(\Phi_{k,u}(x))\leq f(x),\]
for almost all $x$ in $V.$ Then $\Phi_{k}^{*}f_{k}$ tends to $f$
in $L^{1}((V,\nu)\times U).$ In particular, there is a subsequence
$\Phi_{k_{j}}^{*}f_{k_{j}}$such that $f_{k_{j}}(\Phi_{k_{j},u}(x))$
converges to $f(x)$ for almost all pairs $(x,u)$ in $V\times U.$
\end{lem}
\begin{proof}
When $\Phi_{k,u}$ is the identity map (and $V=M)$ the lemma was
essentially obtained in \cite{berm2}. But for completeness we recall
the argument: By the assumption (i)\[
\limsup_{k}\int_{V}\left|f_{k}-f\right|d\nu=2\limsup_{k}\int_{V}\chi_{+,k}(f_{k}-f)d\nu,\]
 where $\chi_{+,k}$ is the characteristic function of the set where
$f_{k}-f$ is non-negative. The right hand side can be estimated by
Fatou's lemma, which (by the dominated convergence theorem) is equivalent
to the inequality \[
\limsup_{k}\int_{X}g_{k}d\nu\leq\int_{X}\limsup_{k}g_{k}d\nu,\]
 if the sequence $g_{k}$ is dominated by an $L^{1}-$function. Taking
$g_{k}=\chi_{k}(f_{k}-f)$ and using the assumption $(ii),$ finishes
the proof of this special case.

Now, for a general map $\Phi_{k,u}$ let $h_{k}$ be the function
on $U\times V$ defined by \[
h_{k}(u,x):=(\Phi_{k,u}^{*}f_{k})(x).\]
The assumption on the Jacobian of $\Phi_{k,u}$ combined with the
assumption $(i)$ shows (by the change of variables formula) that
$(i)$ also holds for $h_{k}$ on $U\times V.$ Indeed, \[
\int_{V}h_{k}d\nu=\int_{V}\Phi_{k,u}^{*}(f_{k}\Phi_{k,u}^{-1*}d\nu)=\int_{\Phi_{k,u}(V)}f_{k}\Phi_{k,u}^{-1*}d\nu\]
and by assumption $\Phi_{k,u}(V)\cap\textrm{supp}f_{k}=V\cap\textrm{supp}f_{k}$
and $\Phi_{k,u}^{-1*}d\nu\rightarrow d\nu$ uniformly, giving \[
\lim_{k}\int_{V}h_{k}d\nu=\lim_{k}\int_{V}f_{k}d\nu=\lim\int_{V}fd\nu\]
In particular, by Fubini's theorem, the corresponding equality holds
over $U\times V$ too. Moreover, by assumption $h_{k}$ is $L^{1}-$dominated.
We may now apply the special case when $\Phi_{k,u}$ is the identity
map to the sequence $h_{k}$ on $U\times V$ and obtain that $h_{k}$
tends to $f$ for almost all pairs $(x,u).$
\end{proof}
In section \ref{sec:Bergman-metric-asymptotics}, we will also have
use for the following

\begin{lem}
\label{lem:int 2}Assume that $(M,\nu)$ is a measure space. Let $f$
and $f_{k}$ be non-negative functions in $L^{1}(M,\nu)$ such that\[
(\textrm{i})\lim_{k}\int_{M}f_{k}d\nu=\int_{M}fd\nu\,\,\,\,\textrm{and}\,\,\,\,(\textrm{ii})\liminf f_{k}(x)\geq f(x)\,\textrm{a.e}\]
Then $f_{k}$ tends to $f$ in $L^{1}(M,\nu).$ In particular, there
is a subsequence $f_{k_{j}}$such that $f_{k_{j}}(x)$ converges to
$f(x)$ for almost all $x.$
\end{lem}
\begin{proof}
Reversing the roles of $f_{k}$ and $f$ in the beginning of the proof
of the previous lemma, the assumption $(i)$ gives \[
\limsup_{k}\int_{X}\left|f_{k}-f\right|d\nu=2\limsup_{k}\int_{X}\chi_{-,k}(f-f_{k})d\nu,\]
where now $\chi_{-,k}$ is the set where $f_{k}<f.$ Let $h_{k}:=\chi_{-,k}(f-f_{k})$
and note that by assumption $(ii)$ we have that \[
\limsup_{k}h_{k}=0\]
 and $h_{k}$ is, by its definition, dominated by the $L^{1}-$function
$f.$ Hence, Fatou's lemma again finishes the proof of the lemma. 
\end{proof}

\section{\label{sec:Bergman-kernel-asymptotics}Bergman kernel asymptotics}

\subsection{Convergence as a current}

\begin{thm}
\label{thm:conv of B}Let $B^{k}$ be the Bergman function for the
Hilbert space $\mathcal{H}_{k}(Y,,L^{k})_{X}$ associated to the polarized
pseudoconcave domain $X$ (section \ref{sub:Notation-and-setup}).
Then \[
k^{-n}B^{k}1_{X}\omega_{n}\rightarrow1_{X}(dd^{c}\phi)_{n}+[\partial X]\wedge\mu\]
as measures on $Y$ in the weak{*}-topology, where $\mu$ is the $2n-1$
form \ref{eq:def of mu} on $\partial X.$
\end{thm}
\begin{proof}
For simplicity we will assume that the restriction of $\rho$ to the
ray close to the boundary where $z$ and the real part of $w$ vanish,
coincides with the restriction of $v$ (the assumption may be removed
as in the proof of proposition 5.5 in \cite{berm3}). Let $B_{X}^{k}(x):=k^{-n}B^{k}(x)$
when $x$ is in the interior region, i.e. $\rho(x)\leq-1/\ln k$ and
$0$ otherwise and let\begin{equation}
B_{\partial X}^{k}(\sigma)=k^{-n}\int_{-\ln k/k}^{0}B^{k}(\sigma,\rho)d\rho\label{eq:pf of conv of B def of B bd}\end{equation}
By formula \ref{eq:ineq int middle} \[
\lim_{k}k^{-n}B^{k}\omega_{n}=\lim_{k}B_{X}^{k}\omega_{n}+[\partial X]\wedge(\lim_{k}B_{\partial X}^{k}(\sigma))d\sigma\]
 as currents. In particular, proposition \ref{pro:Morse eq} gives\[
\lim_{k}\int_{X}B_{X}^{k}\omega_{n}+\lim_{k}\int_{\partial X}B_{\partial X}^{k}(\sigma)d\sigma=\int_{X}(dd^{c}\phi)_{n}+\int_{\partial X}\mu.\]
 Hence, the inequalities \ref{eq:pf of conv of K peak bd} and \ref{eq:ineq B bd}
show that the first and second term on the left hand side in the previous
formula is equal to the first term and the second term, respectively,
in the right hands side. Finally, lemma \ref{lem:integrat theory}
applied to the spaces $X$ and $\partial X$ gives \begin{equation}
\begin{array}{lrcl}
(i) & \lim_{k}B_{X,R}^{k}\omega_{n} & = & (dd^{c}\phi)_{n}\,\textrm{a.e on}\, X\\
(ii) & (\lim_{k}B_{\partial X}^{k}(\sigma)d\sigma) & = & \mu\textrm{ a.e on}\,\partial X.\end{array}\label{eq:pf of conv of B}\end{equation}
This proves the proposition.
\end{proof}
Now we can prove the convergence as a current of the whole Bergman
kernel, stated as theorem \ref{thm:conv of K to current} in the introduction.

\subsubsection*{Proof of theorem \ref{thm:conv of K to current}.}

First observe that \begin{equation}
\lim_{k}k^{-n}\int_{X\times X}f(x,y)\left|K^{k}(x,y)\right|_{k\phi}^{2}\omega_{n}(x)\wedge\omega_{n}(y)=\lim_{k}I_{1,k}+\lim_{k}I_{2,k},\label{eq:pf conv of K I plus I}\end{equation}
where $I_{1,k}$ and $I_{2,k}$ are the integrals obtained by restricting
the integration to the set of all $(x,y)$ such that $\rho(x),\rho(y)\leq-1/\ln k$
and $\rho(x),\rho(y)\geq-\ln k/k^{1},$ respectively. Indeed, if $A_{k}$
denotes the middle region, i.e. the set of all $x$ such that $-1/\ln k\leq\rho(x)\leq-\ln k/k,$
the absolute value of the difference between the integrals in left
hand side and the right hand side in \ref{eq:pf conv of K I plus I}
may be estimated by \[
C\int_{A_{k}}(\int_{X}\left|K^{k}(x,y)\right|_{k\phi}^{2}\omega_{n}(y))\omega_{n}(x)\leq C\int_{A_{k}}B^{k}(x)\omega_{n}(x)\]
using the reproducing property \ref{(I)repr property} of $K_{x}^{k}(y)$
applied to $\alpha_{k}=K_{x}^{k}(y)$ in the last step. By \ref{eq:ineq int middle}
the latter integral tends to zero, when $k$ tends to infinity.

\begin{proof}
Since the proof that \begin{equation}
\lim_{k}I_{1,k}=(\int_{X}f(x,x)(dd^{c}\phi)_{n}\label{eq:pf of conv of K I1}\end{equation}
is completely analogous to the case when $X$ is closed (theorem 2.4
in \cite{berm2}), we will just sketch it here. Take a sequence of
sections $\alpha_{k},$ where $\alpha_{k}$ is a normalized extremal
at the interior point $x.$ Combining the inequality \ref{eq:Intro local weak pointwise}
with the equality $(i)$ in formula \ref{eq:pf of conv of B} shows
that, unless $x$ lies in a set of measure zero, there is a subsequence
of $\alpha_{k}$ such that \begin{equation}
\lim_{k}\left\Vert \alpha_{k}\right\Vert _{k\phi,F_{k}(D_{\ln k})}^{2}=1,\label{(I) claim: peak sections}\end{equation}
restricting the norms to $F_{k}(D_{\ln k}),$ the polydisc of radius
$\textrm{ln\,}k/k$ centered at $x.$ Using the identity \ref{eq:K as extremal}
the integral over $y$ in $I_{1,k,R}$ (for a fixed point $x)$ equals
\[
k^{-n}B^{k}(x)\int_{\rho(y)\leq-1/\ln k}f(x,y)\left|\alpha_{k}(y)\right|_{k\phi}^{2}\omega_{n}(y).\]
Since by \ref{(I) claim: peak sections} the function $\left|\alpha_{k}(y)\right|_{k\phi}^{2}$
converges to the Dirac measure at $x$ in the weak{*}-topology, formula
\ref{(I) claim: peak sections} then proves \ref{eq:pf of conv of K I1}.

Similarly to prove \begin{equation}
\lim_{k}I_{2,k}=(\int_{\partial X}f(x,x)\mu\label{eq:pf of conv of K I2}\end{equation}
 first note that in the limit $f$ may clearly be replaced by its
restriction to $(\partial X)^{2}.$ Replace $x$ in the previous argument
with $x_{k}=(\sigma_{x},v/k)$ and observe that \begin{equation}
\lim_{k}\left\Vert \alpha_{k}\right\Vert _{k\phi F_{k}(\Delta_{\ln k})}^{2}=1,\label{eq:pf of conv of K peak bd}\end{equation}
 restricting the norms to the polydisc of radius $\ln k$ scaled by
the map $F_{k}$ (in formula \ref{eq:scal map bd}). To see this note
that combining $(ii)$ in formula \ref{eq:pf of conv of B} with the
inequality \ref{eq:ineq B bd} gives, using lemma \ref{lem:integrat theory}
applied to the infinite ray $\left\{ \sigma\right\} \times[0,\infty[,$
\[
\lim_{k}k^{-(n+1)}B^{k}(\sigma,v/k)=B^{0}(0,v)\]
for almost all $(\sigma,v).$ The inequality \ref{eq:Morse ineq S bd}
then gives \ref{eq:pf of conv of K peak bd} as in the interior case.
The integral over $y$ in the definition of $I_{2,k}$ (for a fixed
point $x)$ now equals \[
k^{-n}B^{k}(\sigma_{x},\rho_{x})\int_{\rho(y)\geq-\ln k/k^{1}}f(\sigma_{x},\sigma_{y})\left|\alpha_{k}(\sigma_{y},\rho_{y})\right|_{k\phi}^{2}\omega_{n}(y).\]
 By \ref{eq:pf of conv of K peak bd} $\alpha_{k}$ is peaked around
$\sigma_{y}=\sigma_{x}$ showing that \ref{eq:pf of conv of K I2}
is the limit of \[
k^{-n}B^{k}(\sigma_{x},\rho_{x})f(\sigma_{x},\sigma_{x})\int\left|\alpha_{k}(\sigma_{y},\rho_{y})\right|_{k\phi}^{2}\omega_{n}(y).\]
Since the integral of $\left|\alpha_{k}(\sigma_{y},\rho_{y})\right|^{2}$
is equal to one in the limit (by \ref{eq:pf of conv of K peak bd})
formula \ref{eq:pf of conv of B} finally proves \ref{eq:pf of conv of K I1}.
\end{proof}

\subsection{Scaling asymptotics}

In this section scaling asymptotics for the Bergman kernels in the
interior of $X$ and at the boundary of $X$ will be obtained. The
scalings are expressed in terms of the local coordinates introduced
in section \ref{sub:int reg} and \ref{sub:bd reg}, respectively.
In the interior case we will use the notation $B^{(k)}(z)=B^{k}(z/k^{1/2})$
and $K^{(k)}(z,w)=K_{z}^{(k)}(w)=K(z/k^{1/2},w/k^{1/2})$ and similarly
in the boundary case, using the scaling map \ref{eq:scal map bd}
in the latter case. Note that we have suppressed the dependence on
the fixed center (which is the point $x$ in the interiour and the
point $\sigma$ at the boundary).

\begin{lem}
\label{lem:ineq for K}. Let $\phi$ be any smooth Hermitian metric
on the line bundle $L$ over $Y$. Then the scaled Bergman kernels
around each fixed interior point $x$ satisfy\[
(i)\,\limsup_{k}\left\Vert k^{-n}K_{z}^{(k)}\right\Vert _{k\phi_{0}}^{2}\leq\left\Vert K_{z}^{0}\right\Vert _{\phi_{0}}^{2}(=B^{0}(z)),\]
 in terms of the model norms (restricted to a polydisc of radius $\textrm{ln }k$
in the left hand side). Moreover, the left hand side is uniformly
bounded by a constant independent of $z.$ Similarly, for each fixed
boundary point $\sigma$\[
(ii)\,\limsup_{k}\left\Vert k^{-(n+1)}K_{z,w}^{(k)}\right\Vert _{k\phi_{0}}^{2}\leq\left\Vert K_{z,w}^{0}\right\Vert _{\phi_{0}}^{2}(=B^{0}(z,w)),\]
in terms of the model norms (restricted to a polydisc of radius $\textrm{ln }k$
in the left hand side). 
\end{lem}
\begin{proof}
By formula \ref{eq:K as extremal}\begin{equation}
\left\Vert k^{-n}K_{z}^{(k)}\right\Vert _{k\phi{}}^{2}=k^{-n}B^{(k)}(z)(k^{-n}\left\Vert \alpha^{(k)}\right\Vert _{k\phi}^{2}),\label{eq:pf ineq for K}\end{equation}
where $\alpha_{k}$ is an extremal at the point $x_{k}=$$z/k^{1/2}$
with global norm equal to one. Hence, \[
\limsup_{k}\left\Vert k^{-n}K_{z}^{(k)}\right\Vert _{k\phi{}_{0}}^{2}\leq\limsup_{k}k^{-n}B^{(k)}(z)\leq B^{0}(z),\]
where we have used the Morse inequality \ref{eq:Intro local weak pointwise}
in the final step. By the reproducing property of the Bergman kernel
(or the analog of formula \ref{eq:pf ineq for K}) in the model case
this proves $(i).$ The proof of $(ii)$ follows along the same lines,
now using the Morse inequities \ref{eq:ineq B bd}. 
\end{proof}
\begin{lem}
\label{lem:scal conv for B}Let $B^{k}$ be the Bergman function for
the Hilbert space $\mathcal{H}_{k}(Y,L^{k})_{X}$ associated to the
polarized pseudoconcave domain $X$ (section \ref{sub:Notation-and-setup}).
Then $B^{k}$ has a subsequence $B_{j}^{k}$ such that for almost
any point $x$ in the interior of $X$ (i.e. $x\in X-E,$ where $E$
has measure zero) the following scaling asymptotics centered at $x$
hold: 

\begin{equation}
(i)\, k_{j}^{-n}B^{(k_{j})}(z)=B^{0}(z)\label{eq:pf of scal conv for B int}\end{equation}
for almost all $z.$ Similarly, for almost any fixed boundary point
$\sigma$ (i.e. $\sigma a\in\partial X-F,$ where $F$ has measure
zero in $\partial X$)\[
(ii)\,\lim_{k}k_{j}^{-(n+1)}B^{(k_{j})}(z,w)=B^{0}(z,w)\]
for almost all $(z,w).$
\end{lem}
\begin{proof}
First consider the interior case. Taking $\Phi_{k}=Id$ (i.e. the
identity map) in lemma \ref{lem:integrat theory} and using the Morse
inequality \ref{eq:Intro local weak pointwise} applied to the center
(i.e. to $k^{-n}B^{(k)}(0)=k^{-n}B^{k}(x))$ proves \ref{eq:pf of scal conv for B int}
when $z=0,$ i.e. that \begin{equation}
\lim_{k}k^{-n}B^{k}(x)=B_{0}(x)\label{eq:ae conv}\end{equation}
a.e. on $X.$ Now fix a point $x_{0}$ in $X$ and a coordinate neighbourhood
$V$ centered at $x_{0}$ that we identify with a subset of $\C^{n}$.
Let $U$ be a ball of fixed radius centered at the origin in $\C^{n}.$
On $V$ we may write \[
B^{(k)}(A(z)u)=B^{k}(\Phi_{k,u}(z)):=B^{k}(z+A(z)u/k^{1/2}),\]
where $A(z)$ is a matrix-valued function and where the center of
the scaling is $z.$%
\footnote{recall that the definition of $B^{(k)}$ involves a choice of coordinates,
that are orthonormal at $z,$ corresponding to multiplying $A(z)$
by a unitary matrix. But it is clearly enough to obtained the scaling
asymptotics for some choice of $A(z).$ %
} On $U$ the matrix $A(z)$ may even be chosen to depend smoothly
on $z.$ Note that the norm of the Jacobian of $\Phi_{k,u}-Id$ is
bounded by a constant times $1/k^{n/2}.$ Now take a smooth function
$\chi$ supported on $U$ such that $\chi=1$ on some neighbourhood
$V(x_{0})$ of the fixed point $x_{0}.$ Let $f_{k}:=\chi k^{-n}B^{k}$
and $f:=\chi B_{0}$ on $V.$ By \ref{eq:ae conv} (and dominated
convergence) we have \[
\lim_{k}\int_{V}f_{k}\omega_{n}=\int_{V}f\omega_{n}\]
Applying lemma \ref{lem:integrat theory} to $f_{k}$ with $(M,\nu)=(X,\omega_{n})$
and $V$ and $U$ as above and using the Morse inequality \ref{eq:ineq B bd}
proves that $f_{k}$ tends to $f$ for almost all $(x,u)$ in $V\times U.$
In particular, $k^{-n}B^{k}(\Phi_{k,u}(x))$ tends to $B_{0}(x)$
for almost all $(x,u)$ in $V(x_{0})\times U.$ By Fubini's theorem
this means that for almost all fixed $x$ in $V(x_{0})$ we have that
$k^{-n}B^{k}(\Phi_{k,u}(x))$ tends to $B_{0}(x)$ for almost all
$u$ in $U.$ But since $x_{0}$ was arbitrary this proves the interior
case.

To prove the boundary case, first note that as in the proof of theorem
\ref{thm:conv of B}, \[
\lim_{k}\int_{\partial X\times[-\ln k,0]}B^{k}(\sigma,v/k)dvd\sigma)=\int_{\partial X\times[-\ln k,0]}B^{0}(\sigma,v)dvd\sigma).\]
The proof now follows along the lines of the interior case, using
the Morse inequality \ref{eq:ineq B bd} for $B^{k}$.
\end{proof}
We now turn to the proof of the scaling convergence of the Bergman
kernel stated as theorem \ref{thm:scal conv for K} in the introduction.

\subsection*{Proof of theorem \ref{thm:scal conv for K}.}

We first consider the interior case. Fix a point $x$ in $X-E$ where
$E$ is the set of measure zero where the convergence in lemma \ref{lem:scal conv for B}
fails. By the uniform bound in lemma \ref{lem:ineq for K}, the sequence
$k_{j}^{-n}K_{z}^{(k_{j})}$ extended by zero converges (after choosing
a subsequence) weakly to an element $\beta_{z}$ in the model space.
Moreover, since $k_{j}^{-n}K_{z}^{(k_{j})}$ is a holomorphic function
the $L^{2}$- bounds may, using Cauchy estimates, be converted to
$\mathcal{C}^{\infty}-$convergence on any given compact set. In particular,
\[
\lim_{j}k_{j}^{-n}K_{z}^{(k_{j})}(z)=\beta_{z}(z).\]
By lemma \ref{lem:scal conv for B} (which is equivalent to the corresponding
asymptotic identity for $K_{z}^{(k)}(z))$ this means that \begin{equation}
\beta_{z}(z)=K_{z}^{0}(z)\label{eq:pf scal conv of K int A}\end{equation}
under the assumptions of the theorem. The inequality $(i)$ in lemma
\ref{lem:ineq for K} combined with the extremal characterization
\ref{(I)extremal prop of B} of the Bergman function $B_{k}$ then
forces \begin{equation}
\beta_{z}(w)=c_{z}K_{z}^{0}(w)\label{eq:pf of scal conv int B}\end{equation}
for each $w,$ where $c_{z}$ is of unit norm for each fixed $z.$
Combining \ref{eq:pf scal conv of K int A} and \ref{eq:pf of scal conv int B}
when $z=w$ shows that $c_{z}=1.$ All in all we deduce that, for
each fixed $z,$ \[
\lim_{j}k_{j}^{-n}K_{z}^{(k_{j})}=K_{z}^{0}\]
 uniformly on any given compact set. The limit has been established
for a certain subsequence of $k_{j}^{-n}K^{(k_{j})},$ but in fact
it implies point-wise convergence of the sequence $k_{j}^{-n}K^{(k_{j})}$
itself since the limiting function is independent of the subsequence
of $k_{j}^{-n}K^{(k_{j})}.$ Moreover, by the uniform bound $(i)$
in lemma \ref{lem:ineq for K} \[
\left\Vert K^{(k_{j})}-K^{0}\right\Vert \leq C\]
in $L^{2}$ for each fixed compact set in $\C^{n}\times\C^{n}.$ Since
the sequence $(K^{(k_{j})}-K^{0})$ is holomorphic on $\overline{\C^{n}}\times\C^{n}$
Cauchy estimates finally may be used again to convert the $L^{2}$-
convergence on $\C^{n}\times\C^{n}$ to $\mathcal{C}^{\infty}-$convergence
on any fixed compact set of $\C^{n}\times\C^{n}$, proving the theorem
in the interior case.

The boundary case follows along the same lines, now using $(ii)$
in lemma \ref{lem:ineq for K} and \ref{lem:scal conv for B}.

\section{\label{sec:Bergman-metric-asymptotics}Bergman metric asymptotics}

Denote by $\Omega_{k}$ the following global $(1,1)-$current on $Y:$
\[
\Omega_{k}:=dd^{c}(\textrm{ln\,$K^{k}(y,y))(=dd^{c}(\sum_{i}\left|\psi_{i}(y)\right|^{2}))$}\]
where $K^{k}(y,y)$ is the Bergman kernel of the Hilbert space $\mathcal{H}_{k}(Y,L^{k})_{X}$
(with orthonormal basis $(\psi_{i}))$, restricted to the diagonal
and identified with a section of $L\otimes\overline{L}$ over $Y.$
Equivalently, $\Omega_{k}$ is the pull-back of the Fubini-Study metric
$\omega_{FS}$ on $\P^{N}(=\P\mathcal{H}_{k}(Y)_{X})$ (compare section
\ref{sec:Model-examples}) under the Kodaira map \[
Y\rightarrow\P\mathcal{H}_{k}(Y)_{X},\,\,\, y\mapsto(\textrm{$\Psi_{1}(y):\Psi_{2}(y)...:\Psi_{N}(y))$ , }\]
where $\textrm{$(\Psi_{i})$}$ is an orthonormal base for $\mathcal{H}_{k}(Y)_{X},$
i.e $\Omega_{k}$ is the (normalized) curvature of the metric $\textrm{ln\,$K^{k}(y,y)$}$
on $L,$ which is the pull-back of the Fubini-Study metric on the
hyper plane line bundle $\mathcal{O}(1)$ over $\P^{N}(=\P\mathcal{H}_{k}(Y)_{X})$.
We will call $\Omega_{k}$ the \emph{$k$th Bergman metric on $Y$
induced by the polarized domain $X.$}

Now fix a point $\sigma$ in $\partial X$ and recall that $B^{0}$
and $K^{0}$ denote the corresponding model Bergman function and kernel,
respectively, on $\C_{y}^{n}$ defined in section \ref{sub:bd reg}. 

\begin{lem}
\label{lem:mod met}Let $T$ be defined as in formula \ref{eq:def of mu}.
Then\[
dd^{c}(\textrm{ln\,$K^{0}(y,y))$}=td(d^{c}\rho_{0})+dd^{c}\phi_{0}+dt\wedge d^{c}\rho_{0}\]
 where $t=\frac{\partial}{\partial\rho_{0}}\textrm{ln}B^{0}(\rho_{0})$
is strictly increasing, mapping $[-\infty,\infty]$ to $[0,T].$ 
\end{lem}
\begin{proof}
Consider $\psi(\rho_{0})=\textrm{ln}B^{0}(\rho_{0})$ as a function
on $\C^{n}$. Then \[
dd^{c}\psi=d(\frac{\partial\psi}{\partial\rho_{0}}d^{c}\rho_{0})=td(d^{c}\rho_{0})+dt\wedge d^{c}\rho_{0},\]
 where we have used Leibniz rule in the last step and the definition
of $t$ above. To prove the last statement above note that $\psi$
is of the general form \[
\psi(y)=\textrm{ln}\int_{K}e^{\left\langle y,t\right\rangle }d\nu(t),\]
 where $y$ and $t$ are vectors in Euclidean $\R^{N}$ and $d\nu(t)$
is a finite measure supported on a compact set $K.$ Hence, $\psi$
is convex and it follows from well-known convex analysis \cite{gro}
that the gradient of $\psi$ maps $\R^{N}$ bijectively onto the interior
of $K.$ When $N=1$ and $K=[a,b]$ it also follows that $-\infty$
and $\infty$ are mapped to $a$ and $b,$ respectively.
\end{proof}
Next, we will prove theorem \ref{thm:Berg metric} about the convergence
of the Bergman metric.

\subsection*{Proof of theorem \ref{thm:Berg metric}}

Let us first prove the weak convergence of the (normalized) sequence
of Bergman volume forms $(dd^{c}(k^{-1}\textrm{ln\,$K^{k}(y,y))_{n}$ }.$
Since the mass of this sequence of measures is bounded (compare \ref{eq:bound on mass}
below), it is by weak compactness enough to show that any subsequence
has another subsequence that converges to the expected limit. Now
given a first choice of subsequence it has itself a subsequence such
that the scaling convergence in theorem \ref{thm:scal conv for K}
holds. Hence, it will be enough to show weak convergence for the latter
subsequence and to simplify the notation we will assume that it is
indexed by $k$ in the following. 

Let $G^{k}(y):=(dd^{c}(k^{-1}\textrm{ln\,$K^{k}(y,y))_{n}/\omega_{n}$ }$
and\[
G_{X}^{k}(y):=1_{\{\rho\leq-1/\ln k\}}G^{k}(y),\,\,\, G_{\partial X}^{k}(y):=\int_{-\ln k/k}^{\ln k/k}G^{k}(\sigma,\rho)d\rho\]
where $\sigma$ denotes a point in $\partial X$ as in the proof of
theorem \ref{thm:conv of B}. To prove the weak convergence, i.e.
that for any smooth {}``test function'' $f:$\begin{equation}
\lim_{k}\int G^{k}f\omega_{n}=\int_{X(0)}f(dd^{c}\phi)_{n}+\int_{\partial X}f\mu\label{eq:pf met a}\end{equation}
it is clearly enough, by decomposing the previous integral into different
regions, to prove the following 

\begin{claim}
\label{cla:pf of met}The following holds:
\end{claim}
\[
\begin{array}{lrclr}
(a) & G_{X}^{k} & \rightarrow & 1_{X(0)}(dd^{c}\phi)_{n}/\omega_{n} & \textrm{in\,$L^{1}(X,\omega_{n})$}\\
(b) & G_{\partial X}^{k} & \rightarrow & \mu/d\sigma & \textrm{in\,}L^{1}(\partial X,d\sigma)\\
(c) & \int_{-1/\ln k\leq\rho\leq-\ln k/k}G^{k}\omega_{n} & \rightarrow & 0\\
(d) & \int_{R_{k}/k\leq\rho}G^{k}\omega_{n} & \rightarrow & 0\end{array}\]

\begin{proof}
First observe that the following holds: \[
(a'):\,\,\, G_{X}^{k}(x)\rightarrow(1_{X(0)}(dd^{c}\phi)_{n}/\omega_{n})(x)\,\,\textrm{a.e.\, on\,$(X(0),\omega_{n})$}\]
Indeed, consider a fixed point $x$ in the interior region and take
local coordinates $z$ centered and orthonormal at $x.$ Write $z=\zeta/k.$
Then the chain rule gives\begin{equation}
\partial\overline{\partial}(k^{-1}\textrm{ln\,$K^{k})=\sum_{i,j}\frac{\partial^{2}\textrm{ln}K^{(k)}(\zeta))}{\partial\zeta_{i}\partial\bar{\zeta_{j}}}dz_{i}\wedge d\bar{z_{j}}$ }\label{eq:pf of met int}\end{equation}
Note that for each fixed $k$ we may replace $K^{(k)}(\zeta)$ in
the formula above by $k^{-n}K^{(k)}(\zeta),$ since $\partial^{2}(\textrm{ln}k^{-n})=0.$
Now, by the scaling convergence $(i)$ in theorem \ref{thm:scal conv for K}
evaluating \ref{eq:pf of met int} at $0$ shows that at almost any
fixed point $x:$\[
\lim_{k}(\partial\overline{\partial}(k^{-1}\textrm{ln\,$K^{k}))_{n}=(\sum_{i}\lambda_{x,i}dz_{i}\wedge d\bar{z_{i}})_{n}=(\partial\overline{\partial}\phi)_{n}$ . }\]
Next, we will show \[
(b'):\,\,\,\liminf_{k}G_{\partial X}^{k}(\sigma)\geq(\mu/d\sigma)(\sigma)\,\,\textrm{a.e.\, on\,$(\partial X,d\sigma)$}\]
To this end fix a point $\sigma$ in $\partial X-E,$ where $E$ is
the set of measure zero where the scaling convergence $(ii)$ in theorem
\ref{thm:scal conv for K} fails. Take local holomorphic coordinates
$(z,w)$ as in section \ref{sub:bd reg}. Recall that $z$ is in $\C^{n-1}$
and $w=u+iv.$ After a change of variables, $G_{\partial X}^{k}$
may be written in the following way:

\begin{equation}
G_{\partial X}^{k}(\sigma)=\int_{-\ln k}^{\ln k}k^{-1}G^{k}(0,iv'/k)dv'\label{eq:pf of met d}\end{equation}
(we are making the same simplifying assumptions on the fixed ray close
to the boundary as in the beginning of the proof of theorem \ref{thm:conv of B}).
Introduce the scaled coordinates $\zeta=$ $(zk^{1/2},wk).$ Then
the integrand above may be written as

\begin{equation}
k^{-1}G^{k}(0,iv'/k)=\det(\frac{\partial^{2}\textrm{ln}K^{(k)}}{\partial\zeta_{i}\partial\bar{\zeta_{j}}})(0,iv').\label{eq:pf met e}\end{equation}
Indeed, from the definition of $G^{k}$ we have that the left hand
side in the formula above may be written as

\[
\det(\frac{\partial^{2}\textrm{ln}K^{(k)}}{\partial\zeta_{i}\partial\bar{\zeta_{j}}})(\zeta)\cdot(k^{-(n+1)}(d\zeta_{1}\wedge d\bar{\zeta_{1}}\cdots)/(dz_{1}\wedge d\bar{z_{1}}\cdots dw\wedge d\bar{w}))\]
By the definition of the scaled coordinates $\zeta$ the second factor
is equal to one and evaluating the expression at $\zeta=(0,iv')$
then proves \ref{eq:pf met e}.

As in the interior case above, we may now multiply $K^{(k)}$ in \ref{eq:pf met e}
by a factor $k^{-(n+1)}$ and apply theorem \ref{thm:scal conv for K}
($ii)$ combined with lemma \ref{lem:mod met} to obtain

\[
\lim_{k}k^{-1}G^{k}(0,iv'/k)dv'=\det(dd^{c}\phi+tdd^{c}\rho)dt,\]
where $t$ is a function of $v'.$ Fatou's lemma combined with the
change of variables $t=t(v')$ in the integral \ref{eq:pf of met d}
then proves $(b').$ 

Now observe that \ref{eq:pf met a} holds when $f=1.$ Indeed, since
$dd^{c}(k^{-1}\ln K^{k})$ represents the first Chern class of $L$
over $Y$ this follows as in the proof of proposition \ref{pro:Morse eq}.
In particular, splitting the integral gives the following \emph{upper}
bounds on the quantities in the claim \ref{cla:pf of met}: \begin{equation}
\begin{array}{c}
\int_{X(0)}(dd^{c}\phi)_{n}+\int_{\partial X}\mu\geq\\
\limsup_{k}\int G_{X}^{k}\omega_{n}+\liminf_{k}\int_{\partial X}G_{\partial X}^{k}d\sigma+\\
\liminf_{k}(\int_{-1/\ln k_{k}\leq\rho\leq-\ln k_{k}/k}G^{k}\omega_{n}+\int_{\ln k_{k}/k\leq\rho}G^{k}\omega_{n})\end{array}\label{eq:bound on mass}\end{equation}
Moreover, the previous bound clearly also holds whith $\limsup$ in
front of any of the other two integrals (as long as the remaing integrals
have $\liminf$ in front of them) . But then the lower bounds in $(a')$
and $(b')$ above combined with Fatou's lemma force \begin{equation}
\lim_{k}\int G_{X}^{k}\omega_{n}=\int_{X(0)}(dd^{c}\phi)_{n},\,\,\,\lim_{k}\int_{\partial X}G_{\partial X}^{k}d\sigma=\int_{\partial X}\mu.\label{eq:pf met lim}\end{equation}
 Hence, $(a)$ and $(b)$ in the claim \ref{cla:pf of met} follow
from combininig these two limits with $(a')$ and $(b')$ above, using
the integration lemma \ref{lem:int 2}. Finally, combing \ref{eq:bound on mass}
and \ref{eq:pf met lim} proves $(c)$ and $(d).$

The statements $(i)$ and $(ii)$ of the theorem follow directly from
theorem \ref{thm:scal conv for K} combined with lemma \ref{lem:mod met}.
\end{proof}

\section{\label{sec:Distribution-of-random}Equilibrium measures}

In this section we will take $X$ to be \emph{any} given compact set
in $Y$ and $\phi$ any given metric on $L\rightarrow Y$ which is
continuous on $X$ (only the {}``restriction'' of $\phi$ to $X$
will be relevant in the sequal). We will, for simplicity, assume that
$L$ is semi-positive, i.e that it admits some smooth metric with
semi-positive curvature (see \cite{berm4b,b-b} for the setup in the
general case). Comparing with the previous sections we will say that
$X$ (or rather the pair $(X,\phi))$ is a \emph{(semi-)polarized}
\emph{set} if $\phi$ is smooth on $X$ with (semi-)positive curvature
form on all of $X.$

\subsection{Equilibrium metrics}

To a general pair $(X,\phi)$ we may associate the following \emph{equilibrium
metric} on $L\rightarrow Y:$ 

\begin{equation}
\phi_{e}(y)=\sup\left\{ \widetilde{\phi}(y):\,\widetilde{\phi}\in\mathcal{L}_{(X,L)},\,\widetilde{\phi}\leq\phi\,\,\textrm{on$\, X$}\right\} .\label{eq:extem metric}\end{equation}
 where $\mathcal{L}_{(X,L)}$ is the class consisting of all (possibly
singular) metrics on $L$ with positive curvature current. Then the
upper semi-continuous (usc) regularization $\phi_{e}^{*}$ is in $\mathcal{L}_{(X,L)}$
and is locally bounded \cite{g-z}. In particular, the Monge-Ampere
measure $(dd^{c}\phi_{e}^{*})^{n}/n!$ is a well-defined positive
measure by the classical work of Bedford-Taylor \cite{kl}, which
is supported on $X$ and called the \emph{equilibrium measure} associated
to $(X,\phi).$ It was recently introduced in the more general global
setting of quasiplurisubharmonic functions by Guedj-Zeriahi \cite{g-z}.

\subsection{\label{sub:Regularity}Regularity}

In case $(X,\phi)$ is a semi-polarized domain one has that $\phi_{e}=\phi$
on the interiour of $X$ \cite{g-z}, hence the non-trivial contribution
to the equilibrium measure then comes from the boundary of $X.$ The
situation when $X$ is all of $Y,$ but $\phi$ is any (typically
non-positively curved) smooth metric on $L$ is studied in \cite{berm4b}.
In the latter case it follows directly that $\phi_{e}$ is usc, i.e.
$\phi_{e}^{*}=\phi_{e}.$ In the general case the latter property
holds precisely when $\phi_{e}^{*}=\phi_{e}$ on $X$ and we will
then say that $(X,\phi)$ is \emph{regular,} using classical terminology\cite{kl}.
In fact, $\phi_{e}^{*}=\phi_{e}$ on $X$ precisely when $\phi_{e}$
is continuous on all of $Y.$ As we will not consider regularity issues
we refer the interested reader to \cite{b-b} for a recent account,
based on the classical work by Siciak and others. For example, when
$X$ is a domain with smooth boundary $(X,\phi)$ is always regular
(as long as $\phi$ is continuous). 

\begin{rem}
Consider the case when $X$ is a polarized domain in $Y$ with smooth
boundary, so that $\phi_{e}^{*}=\phi_{e},$ which is equal to $\phi$
on $X.$ Then $(dd^{c}\phi_{e})^{n}/n!=0$ on the complement of $X$
\cite{g-z} and by the {}``domination principle'' \cite{b-b} $\phi_{e}$
may then be caracterized as the unique extension of $\phi$ from $X$
to all of $Y$ which solves the Dirichlet problem for the Monge-Ampere
operator on $Y-X.$ Using this it should be possibly to obtain the
convergence of the Monge-Ampere operators in theorem \ref{thm:ln k is equil},
in the special case when $X$ is polarized pseudoconcave domain with
compatible curvatures, from theorem \ref{thm:Berg metric} (compare
the approach in \cite{p-s}).
\end{rem}

\subsection{\label{sub:Bernstein-Markov-measures-and}Bernstein-Markov measures
and general Bergman kernels}

Extending classical terminology (compare \cite{b-l2,b-b}) a measure
$\nu$ is said to satisfy the \emph{Bernstein-Markov property} \emph{with
respect to $(X,\phi)$} if for any positive number $\epsilon$ there
is a constant $C_{\epsilon}$ such that the following inequality holds
for all positive integers $k:$

\begin{equation}
\sup_{x\in X}\left|\alpha_{k}\right|_{k\phi}^{2}(x)\leq C_{\epsilon}e^{k\epsilon}\int_{Y}\left|\alpha_{k}\right|_{k\phi}^{2}\nu\label{eq:bm}\end{equation}
for any element $\alpha_{k}$ of $H^{0}(Y,L^{k}).$ Given such a measure
$\nu$ one obtains a Hilbert space structure on $H^{0}(Y,L^{k})$
by replacing the measure $1_{X}\omega_{n}$ in formula \ref{eq:norm restr}
with $\nu.$ We will denote by $K^{k}$ the corresponding Bergman
kernel, which hence depends on $(\nu,\phi)$ and by $B^{k}\nu$ the
corresponding Bergman measure.

For example, if $X$ is a smooth domain and $\nu=1_{X}\omega_{n},$
where $\omega_{n}$ is a smooth volume form on $Y$ then $\nu$ has
the Bernstein-Markov property with respect to $(X,\phi)$ (compare
\cite{b-b} where this is proved by adapting classical arguments of
Siciak and others).

\subsection{The proof of theorem \ref{thm:ln k is equil}}

In the proof of the theorem \ref{thm:ln k is equil} we will make
use of the following well-known extension lemma, which follows from
the Ohsawa-Takegoshi theorem (compare \cite{bern2}):

\begin{lem}
\label{lem:o t} Let $(L,\phi')$ be a (singular) Hermitian line bundle
such that $\phi'$ has positive curvature form and let $(A,\phi_{A})$
be an ample line bundle with a smooth (but not necesserly positively
curved) metric $\phi_{A}.$ Then the following holds after replacing
$A$ by a sufficiently high tensor power: for any point $y$ where
$\phi'\neq\infty,$ there is an element $\alpha$ in $H^{0}(Y,L^{k}\otimes A)$
such that \begin{equation}
\left|\alpha_{k}(y)\right|_{k\phi'}=1,\,\,\,\left\Vert \alpha_{k}\right\Vert _{Y,k\phi'}\leq C.\label{eq:lem o t}\end{equation}
The constant $C$ is independent of the point $y$ and the power $k.$
\end{lem}
The next lemma is used to reduce the general case to the case when
$L$ is an ample line bundle.

\begin{lem}
\label{lem:decreas}Let $(A,\phi_{A})$ be a positive smooth Hermitian
line bundle and let $\phi_{m}:=m\phi+\phi_{A}$ be the induced metric
on $L^{m}\otimes A$. Then \begin{equation}
((m\phi+\phi_{A})_{e}-\phi_{A})/m\rightarrow\phi_{e}\label{eq:statement of lemma decr}\end{equation}
is a decreasing limit. Moreover, applying the Monge-Ampere operator
to both sides gives a weakly convergent sequence of measures on $Y.$
\end{lem}
\begin{proof}
To simplify the notation we set $\epsilon=1/m$ and interpretate $\phi+\epsilon\phi_{A}$
as a $\Q-$metric on the $\Q$- line bundle $A^{\epsilon}\otimes L$
(these can be defined either in analogy with $\Q-$divisors \cite{la},
or in terms of quasiplurisubharmoic functions as in remark \ref{rem:[notation-in-G-Z]}).
Then $\Phi_{\epsilon}:=(\phi+\epsilon\phi_{A})_{e}-\phi_{A}$ is the
sequence of metrics on $L$ given by the left hand side in \ref{eq:statement of lemma decr}.
To see that $\Phi_{\epsilon}$ decreases as $\epsilon$decreases to
$0$ it is clearly equivalent to prove \[
\epsilon\geq\epsilon'\,\Rightarrow\,(\phi+\epsilon\phi_{A})_{e}\geq(\phi+\epsilon'\phi_{A})_{e}+(\epsilon-\epsilon')\phi_{A}.\]
But this follows since the right hand side is a contender for the
sup defining $(\phi+\epsilon\phi_{A})_{e}$ (using that $dd^{c}\phi_{A}\geq0).$
As a consequence $\lim_{\epsilon\rightarrow0}\Phi_{\epsilon}$ exists
and \[
\lim_{\epsilon\rightarrow0}\Phi_{\epsilon}\geq\Phi_{0}=\phi_{e},\,\,\,\,\, dd^{c}(\lim_{\epsilon\rightarrow0}\Phi_{\epsilon})\geq0\]
Moreover, by definition $\Phi_{\epsilon}\leq(\phi+\epsilon\phi_{A})-\epsilon\phi_{A}\leq\phi$
on $\textrm{int}(X).$ Hence $\lim_{\epsilon\rightarrow0}\Phi_{\epsilon}\leq\phi_{e}$
by the extremal defintion of $\phi_{e}.$ This proves \ref{eq:statement of lemma decr}.
Finally, writing \[
(dd^{c}\Phi_{\epsilon})^{n}/n!=\sum_{k=0}^{n}(-\epsilon)^{(n-k)}(dd^{c}(\phi+\epsilon\phi_{A})_{e})^{k}/k!(dd^{c}\phi_{A})^{(n-k)}/(n-k)!\]
and using that the operator that takes $\psi$ to the current $(dd^{c}\psi)^{k}$)
is continuous \cite{kl} when applied to a decreasing limit of plurisubharmonic
functions (here given by a local representation of $(\phi+\epsilon\phi_{A})_{e}),$
proves the last statement of the lemma. 
\end{proof}
Now we can prove theorem \ref{thm:ln k is equil} stated in the introduction,
saying that the metric on $L$ induced by the Bergman kernel converges
to $\phi_{e}.$

\subsection*{Proof of theorem \ref{thm:ln k is equil}}

Let us first prove the upper bound on $k^{-1}\textrm{ln\,$K^{k}$}$
if $\mu$ is a measure with the Bernstein-Markov property. To this
end fix $\epsilon>0$ and observe that by the very definition of the
latter property and the extremal caractererization \ref{(I)extremal prop of B}
of the Bergman kernel \[
k^{-1}(\textrm{ln\,$K^{k}(x,x)-\textrm{ln\,}C_{\epsilon})-\epsilon\leq\phi$}\]
for $x\in X,$ for any $k.$ In particular,\begin{equation}
k^{-1}(\textrm{ln\,$K^{k}(y,y)-\textrm{ln\,}C_{\epsilon})-\epsilon\leq\phi_{e}$}\label{eq:pf pr ln k up}\end{equation}
on all of $Y,$ by the extremal definition of $\phi_{e},$ which proves
the upper bound corresponding to \ref{eq:statement theorem ln k}. 

Next, let us prove the lower bound in the case $(i).$ Given $\epsilon>0$
fix an arbitrary point $y$ in $Y.$ A standard regularization argument
involving Demailly's regularization theorem (see \cite{b-b}) yields
a candidate $\phi'$ for the sup defining $\phi_{e}$ such that $\phi'$
is continuous on $X$ and $\phi'(y)\geq\phi_{e}(y)-\epsilon.$ Now
take an element $\alpha_{k}$ in $H^{0}(Y,L^{k}\otimes A)$ furnished
by lemma \ref{lem:o t}. Since by construction $\phi'\leq\phi$ on
$X$ we have that \begin{equation}
\int_{X}\left|\alpha_{k}\right|_{k\phi+\phi_{A}}^{2}\omega_{n}\leq\int_{Y}\left|\alpha_{k}\right|_{k\phi'+\phi_{A}}^{2}\omega_{n}\label{eq:pf of pr norm gl}\end{equation}
Hence, \ref{eq:lem o t} in lemma \ref{lem:o t} gives that \[
k^{-1}\textrm{ln\,($\left|\alpha_{k}(y)\right|_{k\phi'+\phi_{A}}^{2}/\int_{X}\left|\alpha_{k}\right|_{k\phi+\phi_{A}}^{2}\omega_{n})\geq C'k^{-1}$ }\]
and since $\phi'(y)\geq\phi_{e}(y)-\epsilon$ we then obtain \[
k^{-1}\textrm{ln\,$(\left|\alpha_{k}(y)\right|^{2}/\int_{X}\left|\alpha_{k}\right|_{k\phi+\phi_{A}}^{2}\omega_{n})\geq\phi_{e}(y)-\epsilon+(\phi_{A}(y)+C')k^{-1}$ }\]
Finally setting $(A,\phi_{A})=(L^{k_{0}},k_{0}\phi)$ for $k_{0}$
a fixed large natural number and writiting $L^{k}=L^{k-k_{0}}\otimes A$
proves the lower bound corresponding to \ref{eq:statement theorem ln k}
in the case $(i).$

The case $(ii)$ now follows from the previous lower bound and the
fact that the Bernstein-Markov inequality used to get the upper bound
may be replaced by the following stronger inequality if $X$ is a
smooth pseudoconcave domain: \begin{equation}
\sup_{x\in X}\left|\alpha_{k}(x)\right|_{k\phi}^{2}/\int_{X}\left|\alpha_{k}\right|_{k\phi}^{2}\omega_{n}\leq Ck^{(n+1)}\label{eq:pf pr ln k sup}\end{equation}
uniformly in $k.$ Indeed, the previous bound is a direct consequence
of the Morse inequalities in section \ref{sec:Morse-(in)equalities-and}
(see also \cite{berm3} for the {}``middle region''). 

To prove the lower bound in the case $(iii),$ fix $\epsilon>0$ and
denote by $X_{\delta}$ the closure of an open $\delta-$neighbourhood
of $X.$ By the previous argument used to prove the lower bound in
the case $(ii)$ applied to $X_{\delta}$ (for $\delta$ sufficently
small) and with $\phi'=\phi_{e}$ (which is usc by assumption) it
is enough to prove that 

\begin{equation}
\int_{X}\left|\alpha_{k}\right|_{k\phi_{e}}^{2}\nu\leq C_{\epsilon}e^{k\epsilon}\int_{X_{\delta}}\left|\alpha_{k}\right|_{k\phi_{e}}^{2}\omega_{n}.\label{eq:pf of thm ln k subm}\end{equation}
But this is a simple consequence of the submean property of holomorphic
functions. Indeed, for any fixed $x$ in $X$ the latter property
gives\[
\left|\alpha_{k}\right|^{2}(x)\leq C_{\delta}\int_{B_{\delta}(x)}\left|\alpha_{k}\right|^{2}\omega_{n},\]
 in a fixed trivialization of $L$ on the coordinate ball $B_{\delta}(x)$
of radius $\delta$ centered at $x.$ Now since, by assumption, $\phi_{e}$
is usc we may chose $\delta$ sufficiently small that \begin{equation}
\left|\alpha_{k}\right|_{k\phi_{e}}^{2}(x)\leq C_{\delta(\epsilon)}e^{k\epsilon}\int_{B_{\delta}(x)}\left|\alpha_{k}\right|_{k\phi_{e}}^{2}\omega_{n}(\leq C_{\delta(\epsilon)}e^{k\epsilon}\int_{X_{\delta}}\left|\alpha_{k}\right|_{k\phi_{e}}^{2}\omega_{n})\label{eq:subm}\end{equation}
Hence, integrating over $x$ proves \ref{eq:pf of thm ln k subm}.

Finally, to prove the convergence \ref{eq:bergman vol as equi} in
the theorem recall that the Monge-Ampere operator (mapping $\phi$
to $(dd^{c}\phi)^{n}$) is continuous when applied to a uniform limit
\cite{kl} of plurisubharmonic functions.

\begin{rem}
\label{rem:l infty}Given a compact subset $X$ and a continuous metric
$\phi$ on an ample line bundle $L,$ let \[
\psi_{k}(y):=\textrm{ln\,$($}\sup_{\alpha_{k}\in H^{0}(Y,L^{k}).}\frac{\left|\alpha_{k}(y)\right|^{2}}{\sup_{Y}\left|\alpha_{k}(y)\right|_{k\phi}^{2}}),\]
 which is an $L^{\infty}$-version of the $k$th Bergman metric. By
definition $\psi_{k}\leq\phi_{e}$ on $X$ for any $k.$ Hence, the
proof of the lower bound in Theorem \ref{thm:ln k is equil} $(iii)$
in gives that $\psi_{k}$ also converges uniformly to $\phi_{e}$
if $(X,\phi)$ is regular. Moreover, a slight modification of the
argument gives that if $X$ is \emph{any} compact set, then the corresponding
\emph{point-wise} convergence always holds. Indeed, given a fixed
point $y$ in $Y$ and $\epsilon>0$ one applies lemma \ref{lem:o t}
to a continuous metric $\phi'$ defined as in the proof of $(i)$
in Theorem \ref{thm:ln k is equil}. Then the convergence will in
general only be point-wise since the $\delta$ in the estimate \ref{eq:subm}
will depend on the {}``oscillation'' of $\phi'$ (which plays the
role of $\phi_{e})$ and hence on the point $y.$ 
\end{rem}

\subsection*{Proof of corollary \ref{cor:expl equil}}

First assume that $L$ is ample. Then combining theorem \ref{thm:Berg metric}
and theorem \ref{thm:ln k is equil} immediately proves the corollary.
When $L$ is not ample we replace it by the ample line bundle $L^{m}\otimes A.$
Finally, letting $m$ tend to infinity and combining the ample case
with lemma \ref{lem:decreas} then finishes the proof of the corollary.

\begin{example}
Under the assumption 2 in section \ref{sub:Assumptions-on-the} ,
we have that $\phi_{e}\equiv0$ on $Y-X$ (also assuming that $L$
is a positive line bundle). Indeed, in this case $\phi_{e}$ is continuous
on all of $Y$ (since it is the uniform limit of continuous functions
according to theorem \ref{thm:ln k is equil}). Hence, $\phi_{e}\equiv0$
on $Y-X$ by the uniqueness of solutions to corresponding Dirichlet
type problem. 
\end{example}
In particular, the previous example shows that in the case considered
in example \ref{ex: lnplus} $\phi_{e}$ corresponds to the function
$\ln(\left|\zeta\right|_{+}^{2})$ in $\C^{n}.$ %
\footnote{This is the \emph{Siciak-Zaharjuta extremal function} of the unit-ball
(also called the \emph{pluricomplex Green function with a pole at
infinity} \cite{kl}). %
} This also follows from the following class of examples:

\begin{example}
In the case considered in example \ref{exa:pull-back bd} we have
that \[
\phi_{e}(z,w):=\phi_{G}(z)+\rho_{+}\]
 on $Y-X,$ using the notation $f_{+}=f$ when $f\geq0$ and $f_{+}=0$
otherwise. In fact, this is the extension of $\phi$ described in
that example and its regularization also described there is incidentally
very similar to the one obtained from the limit of $k^{-1}\textrm{(ln\,$K^{k}(y,y)$ }$(which
may be computed following the constant curvature case in section 4
in \cite{berm3}). Note that $\phi_{e}$ is continuous up to the boundary
of $Y-X$ and solves the Dirichlet type problem in formula \ref{eq:propert of eq}. 
\end{example}
\begin{rem}
\label{rem:[notation-in-G-Z]}In the more general setup of Guedj-Zeriahi
\cite{g-z} the pseudoconcave domain $X$ is replaced by any given
Borel set $K$ in $Y.$ Moreover, the function $V_{K,\omega}:=\phi_{e}(y)-\phi(y)$,
where $\omega=dd^{c}\phi,$ is called the \emph{global extremal function}
in \cite{g-z} and it is an example of an \emph{$dd^{c}\phi-$plurisubharmonic}
function. In particular, the (generalized) equilibrium measure introduced
above may be written as $(dd^{c}\phi_{e})_{n}=(\omega+dd^{c}V_{K,\omega})_{n}$. 
\end{rem}

\section{\label{sec:Open-problems}Open problems}

In view of Theorem \ref{thm:conv of B}, concerning a polarized domain
$X$ with compatible curvatures, it is natural to make the following
conjecture : 

\begin{conjecture}
\label{con:a}Suppose that $(X,\phi)$ is regular and that the measure
$\nu$ \emph{}has the Bernstein-Markov property \emph{}with respect
to $(X,\phi)$. \emph{}Then the normalized Bergman measure $k^{-n}B^{k}\nu$
converges weakly to the the equilibrium measure $(dd^{c}\phi_{e})_{n}$
associated to $(X,\phi).$
\end{conjecture}
If the Bergman measure in the conjecture above is replaced by the
normalized volume form $(dd^{c}(k^{-1}\textrm{ln\,$K^{k}(y,y))_{n}$ }$
of the corresponding Bergman metric (compare section \ref{sec:Bergman-metric-asymptotics})
then the statement does hold, according to theorem \ref{thm:ln k is equil}.
Hence, the conjecture above is equivalent to the following

\begin{conjecture}
\label{con:b}Under the assumptions of the previous conjecture the
weak limits of $k^{-n}B^{k}(y)\nu$ and $(dd^{c}(k^{-1}\textrm{ln\,$K^{k}(y,y))_{n}$ }$
both exist and coincide.
\end{conjecture}
In the {}``weighted classical setting'' in $\C^{n}$ the conjecture
\ref{con:a} was made independently very recently by Bloom-Levenberg
in \cite{b-l3}. In \cite{b-l2} Bloom-Levenberg proved the conjecture
for $n=1,$ i.e. in the complex plane, by using the fact that in this
case the equilibrium measure may be characterized as a minimizer of
the weighted logarithmic energy (see \cite{b-l2} for further references
concerning the planar case). When $X=Y,$ the metric $\phi$ is a
smooth metric and $\mu$ is a smooth volume form the conjecture was
proved very recently in \cite{b-l2} (with $L$ \emph{any} line bundle
over $X).$%
\footnote{In fact, a much stronger convergence result was obtained in this case,
by showing that the corresponding equilibrium measure is absolutely
continuous with respect to the Lebesgue measure and is equal to the
limit of $k^{-n}B^{k}$ almost everywhere on $X.$%
} By reducing to this case the conjecture can be shown to hold when
$X$ is any disc subbundle of a $\P^{1}$-bundle (as in section \ref{sec:Model-examples},
but without any curvature assumptions). The proof will be given in
\cite{berm5}.

Returning to the case when $X$ is a polarized pseudoconcave domain
and $\nu$ a volume form note that if conjecture \ref{con:a} holds
than one obtains the following bound on the corresponding equilibrium
measure from the local Morse inequalities in section \ref{sec:Morse-(in)equalities-and}:
\[
(dd^{c}\phi_{e})_{n}\leq1_{X}(dd^{c}\phi)_{n}+[\partial X]\wedge\mu\]
Finally, it seems also natural to conjecture that the factors $C_{\epsilon}e^{k\epsilon}$in
\ref{eq:bm} may be replaced by $Ck^{n+1}$ for some constant $C$
if $X$ is smooth domain in $Y$ and $\phi$ is a smooth metric on
$L.$ See \cite{zer} for results in this direction in the {}``classical
setting''. As shown in the following appendix the latter conjecture
does hold when $X$ is pseudoconcave. One can also ask if the equilibrium
measure of a polarized domain with smooth boundary is such that the
measure $(dd^{c}\phi_{e})_{n}-1_{X}(dd^{c}\phi)_{n}$ which is supported
on $\partial X$ is absolutely continuous w.r.t the {}``surface measure''
on $\partial X$?

\section{Appendix }

\subsection*{\label{app: A}Boundary estimates}

In this section we will give the proof of some boundary estimates
for the scaled Bergman function $B^{(k)},$ refered to in section
\ref{sub:bd reg}. The arguments are essentially contained in \cite{berm3},
but for completeness we provide some elemantary arguments, which don't
use any subelliptic estimates (as opposed to \cite{berm3}). The notation
in section \ref{sub:bd reg} will be used, but for notational convenience
we assume that the curvature eigenvalues $\mu_{i}$ are all equal
to $-1.$ Moreover, note that on $X_{k}$ the weighted norm is equivalent
to the unweighted one: for any function $f$ we have\[
(1/C)\left\Vert f\right\Vert _{X_{k}}^{2}\leq\left\Vert f\right\Vert _{X_{k},k\phi}^{2}\leq C\left\Vert f\right\Vert _{X_{k}}^{2},\]
 which follows immediately from the convergence of the scaled metrics.
The following lemma uses the pseudoconcavity of $X_{0}$ to estimate
the values of a holomorphic function $f$ on a polydisc centered at
$0,$ with the norm of $f$ inside $X_{k}.$ 

\begin{lem}
Let $f$ be a holomorphic function on $X_{k}.$ Then \[
\sup_{\Delta_{R}}\left|f\right|^{2}\leq C_{R}\left\Vert f\right\Vert _{X_{k}}^{2}\]
 
\end{lem}
\begin{proof}
By a scaling argument we may asume that $R=1.$ First observe that,
by the submean property of holomorphic functions, \[
\sup_{\Delta}\left|f\right|^{2}\leq C\left\Vert f\right\Vert _{\Delta}^{2}.\]
 Hence, it is enough to prove that the integral outside $X_{k}$ may
be estimated by an integral inside of $X_{k}:$\begin{equation}
\left\Vert f\right\Vert _{\Delta-X_{k}}^{2}\leq C'\left\Vert f\right\Vert _{X_{k}}^{2}.\label{eq:app mass outside}\end{equation}
The latter estimate is essentially a well-known consequence of the
pseudoconcavity of $X.$ To see this, note that for any given point
$(z_{0},w_{0})$ in $\Delta$ we have that \begin{equation}
S_{0}:=\left\{ (z,w_{0}):3>\left|z-z_{0}\right|>2\right\} \subseteq X_{k}\label{eq:app s0}\end{equation}
 for $k$ large. Indeed, by the uniform convergence \ref{eq:conv of domains}
it is enough to prove the inclusion into $X_{0},$ which in turn follows
from the bound\[
\rho_{0}(z,w_{0}):=\textrm{Im}w_{0}-\left|z\right|^{2}\leq1-\left|z\right|^{2}<0,\]
 if $(z,w_{0})$ is in $S_{0}.$ Next note that by the submean-property
of the holomorphic function $f(\cdot,w_{0})$ we have\[
\left|f\right|^{2}(z_{0},w_{0})\leq C\frac{\int_{2}^{3}(\int_{\left|z\right|=r}\left|f\right|^{2}(z,w_{0})d\widehat{\sigma})r^{2n-1}dr}{\int_{2}^{3}r^{2n-1}dr}=C'\int_{S_{0}}\left|f\right|^{2}(z,w_{0})dz\wedge d\bar{z},\]
 where $d\widehat{\sigma}$ is the normalized measure on the sphere
in $\C_{z}^{n-1}$ of radius $r$ centered at $z_{0}.$ Finally the
bound \ref{eq:app mass outside} is obained by first integrating over
the $z-$variable in the left hand side of \ref{eq:app mass outside}
and then using the previous point-wise estimate on the integrand.
The point is that, by \ref{eq:app s0}, $S_{0}$ is a subset of $X_{0}.$
\end{proof}
The next lemma is independent of any curvature assumptions.

\begin{lem}
Let $f$ be a holomorphic function on $X_{k}.$ Then for $v\in[-\frac{1}{2}\ln k,0]$
we have that \[
\left|f(0,iv)\right|^{2}\leq\frac{C}{v^{2}}\left\Vert f\right\Vert _{X_{k}}^{2}\]

\end{lem}
\begin{proof}
By the submean propery of the holomorphic function $f(0,\cdot)$ we
have \[
\left|f(0,iv)\right|^{2}\leq\frac{C}{v^{2}}\int_{\left|w-iv\right|\leq-v/2}\left|f\right|^{2}(0,w)dw\wedge d\bar{w}.\]
 Note that, by assumption, the integration takes place over points
inside $X_{k}.$ Finally, estimating $\left|f\right|^{2}(0,w),$ using
the submean property of $f(\cdot,w)$ over the unit-ball in the $z-$variables,
then finishes the proof of the lemma.
\end{proof}
Now by combining the two previous lemmas we obtain that the function
$1_{R_{k}\leq v\leq0}(v)B^{(k)}(0,iv)$ is dominated by the following
function which is in $L^{1}]-\infty,0[:$ \[
g(v):=C,\,\textrm{for\,}v\geq-1\,\,\, g(v):=\frac{C}{v^{2}}\,\textrm{for\,}v<-1\]

\end{document}